\documentclass[a4paper, 11pt]{amsart}
\usepackage{amssymb}
\usepackage{amsmath,amscd}
\usepackage{amsthm}
\usepackage{mathtools}
\usepackage[mathscr]{euscript}
\usepackage[all]{xy}
\usepackage{color}
\usepackage[dvipsnames]{xcolor}
\usepackage[utf8]{inputenc}
\normalfont
\usepackage[T1]{fontenc}
\usepackage[textwidth=14cm,hcentering]{geometry}

\usepackage[colorlinks=true,linkcolor=blue,citecolor=red]{hyperref}
\usepackage{enumerate}
\usepackage[shortlabels]{enumitem}
\usepackage{stackrel}
\usepackage{comment}

\usepackage{mdwlist}
\usepackage{color}
\usepackage[dvipsnames]{xcolor}
\usepackage{array} 
\input xy
\xyoption{all}

\usepackage{tikz}
\usepackage{tikz-cd}
\usetikzlibrary{arrows,calc,matrix,topaths,positioning,scopes,shapes,decorations}
\usepackage{graphicx}
\newcommand{\kkpic}{\vcenter{\hbox{\tiny$\bullet$}}} 
\newcommand{\kkpnpic}{*} 

\usepackage{enumitem}
\usepackage{ stmaryrd }
\usepackage{euscript}
\usepackage{ upgreek }

\makeatletter
\@namedef{subjclassname@2020}{\textup{2020} Mathematics Subject Classification}
\makeatother

\newtheorem{prop}{Proposition}[section]
\newtheorem{teo}[prop]{Theorem}
\newtheorem{lem}[prop]{Lemma}
\newtheorem{cor}[prop]{Corollary}

\theoremstyle{definition}
\newtheorem{defi}[prop]{Definition}
\newtheorem{example}[prop]{Example}

\newtheorem{rmk}[prop]{Remark}
\newtheorem{notation}[prop]{Notation}

\numberwithin{prop}{subsection} 

\newcommand{\ZZ}{\mathbb{Z}}

\newcommand{\Dd}{\mathcal{D}}
\newcommand{\Ee}{\mathcal{E}}
\newcommand{\Ff}{\mathcal{F}} 

\newcommand{\Kk}{\mathcal{K}}

\newcommand{\Nn}{\mathcal{N}} 

\newcommand{\Rr}{\mathcal{R}}
\newcommand{\Qq}{\mathcal{Q}} 
\newcommand{\Ss}{\mathcal{S}}
\newcommand{\Sf}{\mathfrak{S}}
\newcommand{\Tt}{\mathcal{T}}
\newcommand{\Uu}{\mathcal{U}}

\newcommand{\Ww}{\mathcal{W}}

\newcommand{\Yy}{\mathcal{Y}}
\newcommand{\Zz}{\mathcal{Z}}

\newcommand{\Fib}{\mathcal{F}ib}
\newcommand{\Cof}{\mathcal{C}\mkern-1mu o\mkern-2mu f}
\newcommand{\Iso}{\mathcal{I}so}

\newcommand{\Hom}{\mathrm{Hom}}
\newcommand{\Id}{\mathrm{Id}}

\newcommand{\Ker}{\operatorname{Ker}}
\newcommand{\Coker}{\operatorname{Coker}}
\newcommand{\Img}{\operatorname{Im}}
\newcommand{\RelCat}{\mathrm{RelCat}}

\newcommand{\Dec}{\operatorname{Dec}}
\newcommand{\LDec}{\operatorname{LDec}}
\newcommand{\Shift}{\operatorname{Shift}}
\newcommand{\colim}{\mathrm{colim}}
\newcommand{\kk}{\mathtt{R}}

\newcommand{\cone}{\operatorname{Cone}}
\newcommand{\tr}{\operatorname{tr}}

\newcommand{\rbc}{r\text{-}\mathsf{C}}            
\newcommand{\ibc}{i\text{-}\mathsf{C}}  
\newcommand{\mbc}{m\text{-}\mathsf{C}_\kk}            
\newcommand{\bimod}{\mathsf{bgMod}_\kk} 		
\newcommand{\modu}{\mathsf{Mod}_\kk} 		
\newcommand{\spse}{\mathsf{SpSe}} 			
\newcommand{\espse}{\mathsf{ESpSe}} 			
\newcommand{\delrcat}{\langle\delta_r\rangle} 
\newcommand{\D}{\mathsf{D}} 
\newcommand{\lwb}{\mathsf{LWB}} 
\newcommand{\lwbe}{\mathsf{(LWB)^e}} 
\newcommand{\lwbs}{\mathsf{(LWB)^s}} 
\newcommand{\Bcat}{\mathsf{B}}
\newcommand{\Ccat}{\mathsf{C}}
\newcommand{\fcpx}{\mathbf{F}\mathrm{C}_\kk} 	

\newcommand{\pb}{\ar@{}[dr]|{\mbox{\LARGE{$\lrcorner$}}}}
\newcommand{\lra}{\longrightarrow}

\newcommand{\cl}[2]{{#1}\mathrm{\hbox{-}{#2}}}

\newcommand{\veq}{\mathrel{\rotatebox{90}{$=$}}}

\newcommand{\vmat}[2]{\tiny{\begin{pmatrix} #1\\ #2\end{pmatrix}}}

\title{Spectral sequences via linear presheaves}

\author{Muriel Livernet}
\address[M. Livernet]{
Université Paris Cité and Sorbonne Université, CNRS, IMJ-PRG, ENS-PSL, DMA F-75013 Paris, France}
\email{muriel.livernet@imj-prg.fr}

\author{Sarah Whitehouse}
\address[S. Whitehouse]{
School of Mathematical and Physical Sciences\\ 
University of Sheffield\\ S3 7RH\\ England}
\email{s.whitehouse@sheffield.ac.uk }

\thanks{}



\begin{document}

\date{\today}

\begin{abstract}
We study homotopy theory of the category of spectral sequences with
respect to the class of weak equivalences given by 
maps which are quasi-isomorphisms on a fixed page.
 We introduce
the category of extended spectral sequences and show that this is bicomplete
by analysis of a certain linear presheaf category modelled on discs. We endow
the category of extended spectral sequences with various model category
structures, restricting to give the almost Brown category structures on spectral sequences
of our earlier work.
 One of these has the property that spectral sequences is a
homotopically full subcategory. By results of Meier, this exhibits
the category of spectral sequences as a fibrant object in the Barwick-Kan 
model structure on relative categories, that is, it gives a model for an infinity category
of spectral sequences. We also use the presheaf 
approach to define two décalage functors on spectral 
sequences, left and right adjoint to a shift functor,
thereby clarifying prior use of the term d\'ecalage in connection with spectral sequences.
\end{abstract}

\maketitle

\setcounter{tocdepth}{2}
\setcounter{secnumdepth}{2}

\tableofcontents

\section{Introduction}

We study the category $\spse$ of spectral sequences and its homotopy theory.  Since $\spse$ is neither complete nor cocomplete,
it does not admit model category structures. In~\cite[Theorem 5.3.1]{LW24}, we established a weaker homotopical framework, that of an \emph{almost
Brown category}, and exhibited such a structure, $\spse_r$,
on $\spse$ for each $r\geq 0$. The
class of weak equivalences is given by 
maps of spectral sequences which are quasi-isomorphisms on page $r$.
Here, we situate $\spse$ as a subcategory of the category $\espse$ of
\emph{extended spectral sequences} and exhibit various model category structures on this category.
This setting provides a new perspective on the category of spectral sequences and its homotopy theory
and we deduce consequences for the infinity category of spectral sequences.

To this end, we introduce and study a category of linear presheaves closely related to the category of spectral sequences $\spse$. This is the category $\lwb$
of \emph{linear witness books}, a linear presheaf category built from suitable disc objects. Intermediate between these two categories is
the category of extended spectral sequences $\espse$. In an extended spectral sequence, we require a specified morphism from each page to the homology of the previous one, but we drop the requirement, crucial to spectral sequences, that this should be an isomorphism.
\medskip

The first parts of the paper are categorical. They are largely motivated by the wish to view $\spse$ as a subcategory of a convenient
bicomplete category. This is the role played by the category of extended spectral sequences $\espse$. There is a choice involved here,
as we weaken the requirement in spectral sequences that a page is isomorphic to the homology of the previous one. Here we simply 
require a map, but not that it be an isomorphism. We have chosen that the maps go from a page to the homology of the previous one.
In some sense, this choice effectively gives preference to colimits over limits. 
It fits well with the notion of \emph{witness} cycles and boundaries appearing in our previous work~\cite{CELW20}.
Another motivation for our choice is that it is likely to be
better behaved in terms of monoidal structure. Although we do not pursue that direction here, we note that Brotherston has shown that, 
on the category of filtered complexes,
model structures which are closely related to spectral sequences are monoidal~\cite{Bro24}. 

It turns out that $\espse$ is bicomplete; see Theorem~\ref{T:espse-lims-colims}.
Colimits are calculated pagewise, but to understand limits is less straightforward. Here the linear presheaf category of linear witness books
plays a vital role, together with a pair of  adjoint functors $(\Qq, \Nn)$. The terminology of linear witness books is chosen because
the objects of this category can be viewed as having pages like those of a spectral sequence, with \emph{witness} maps from a page
to the previous one, as well as degeneracy maps in the other direction. This means objects have extra data, compared with spectral sequences,
witnessing how elements end up on the $r$-page.

We establish an adjunction $\Qq\dashv \Nn$ of functors between $\lwb$ and $\espse$ and use its properties to identify subcategories $\lwbe$ and $\lwbs$ of $\lwb$ 
equivalent to $\espse$ and $\spse$ respectively. As $\lwbe$ is a full reflective subcategory of $\lwb$ 
it has all (small) limits and colimits and thus so does $\espse$.

This setting offers insight into \emph{d\'ecalage} for spectral sequences. We study 
truncation functors on the underlying category $\D$ on which we take our linear presheaves.
The embedding of a suitably truncated version of $\D$ into $\D$
has both a left and a right adjoint. This triple of adjoint functors gives rise to a chain of five adjoint functors on $\lwb$. Of these, we
 will see that the leftmost three are internal to $\lwbe$ and 
 so there is a corresponding triple of adjoints
  on $\espse$. These also restrict to $\spse$. 
  In particular we obtain a shift functor on $\espse$ or $\spse$ with both a left adjoint $\LDec$ and a right adjoint $\Dec$,
  two versions of d\'ecalage; see Theorem~\ref{T:ldec-shift-dec}.
  The terminology is explained by noting that these functors are suitably compatible with Deligne's functors $\Dec^*$ and $\Dec$ for filtered 
complexes~\cite{DeHII}.  
D\'ecalage has been studied for more general filtered objects, notably in recent work of Antieau~\cite{antieau2024} which focuses on the role of Beilinson \(t\)-structures; see also~\cite{Hed20}. Although d\'ecalage in this sense is always closely connected to the study of spectral sequences, we are not aware of
other work defining d\'ecalage functors directly on the category of spectral sequences, as we do here. Nonetheless, reference is quite often made to d\'ecalage of
a spectral sequence and there are important instances of this relationship. For example, Rognes notes in~\cite{Rog} that it is common to call the Whitehead tower spectral
sequence the d\'ecalage of the Atiyah–Hirzebruch spectral sequence; Levine establishes such a relationship between the Adams-Novikov spectral sequence and Voevodsky's slice tower spectral sequence~\cite{Lev15} and Burklund-Hahn-Senger do so for the
$C_2$-effective slice spectral sequence and the $MU_\mathbb{R}$-based Adams–Novikov spectral sequence~\cite{BHS}.

\medskip
Later parts of the paper are homotopical, extending the study of the homotopy theory of spectral sequences initiated in
~\cite{LW24}. For each $r\geq 0$, we study spectral sequences with $r$-quasi-isomorphisms as a \emph{relative category}, denoted $(\spse, \Ee_r)$.
Having situated spectral sequences inside the bicomplete category of extended spectral sequences $\espse$, we now establish model category structures there, restricting to the relevant structure on spectral sequences. 

For each $r\geq 0$, we show the existence of a model category structure on $\espse$, restricting to
recover the corresponding underlying almost Brown category structure on $\spse$. Indeed there are two flavours of
such model structures; in both the fibrations are those maps $f$ that are surjective on pages $0$ to $r$. Theorem~\ref{T:modelespse2}
gives a structure $\espse_r$ where the weak equivalences are those maps such that the component on the $r$-page is a quasi-isomorphism. In
Theorem~\ref{T:modelespse3}, we obtain $\espse'_r$, where the weak equivalences are those maps such that the component on the 
$r$-page is a quasi-isomorphism and the components on all higher pages are isomorphisms. As relative categories, both $\espse_r$ and $\espse'_r$ restrict to 
$(\spse, \Ee_r)$. 

The methods of proof use the category $\lwb$.  
For the first family of structures, we obtain a cofibrantly generated model category structure on $\lwb$ by transfer of a projective-type 
model structure on the category of $r$-bigraded complexes and then modify this in order to produce a version $\lwb_r$ which 
is closely related to the relevant structure on $\spse$. This model structure is then transferred to produce $\espse_r$.
 The second family of model structures is established by directly checking the axioms, 
making use of the existence of the first family. 

These model structures have the following relationships to each other.
We show in Corollary~\ref{C:w-dec-espse} that the model categories $\espse_r$ for different $r$ are all Quillen equivalent via shift and d\'ecalage functors
 and indeed, these
are all Quillen equivalent to a projective-type model category structure on the category of $0$-bigraded complexes. Similarly,
 the model categories $\espse'_r$ for different $r$ are all Quillen equivalent, see Proposition~\ref{P:w-dec-espse-new}. 
Proposition~\ref{P:localization}
shows that the identity functor $\espse'_r\to \espse_r$ is a right Bousfield localization which is not a Quillen equivalence.
Thus we provide a right \emph{delocalization}
of (a model category Quillen equivalent to) the projective model structure on $0$-bigraded complexes, a bigraded version of chain complexes.

The model category $\espse'_0$ has the important feature that $\spse_0$ is
a homotopically full subcategory, in the sense that any object weakly equivalent to a spectral sequence
is itself a spectral sequence. 
In~\cite{BK} Barwick and Kan provide a model category structure on the category of relative categories, Quillen equivalent to
Rezk's complete Segal space model structure on simplicial spaces~\cite{Rezk}, thus establishing
another model for a homotopy theory of homotopy theories. 
 Results of Meier~\cite[Theorem 4.13]{Meier16} allow us to conclude, in Theorem~\ref{T:fib-recat}, that 
$(\spse, \Ee_0)$ is a fibrant relative category in this model.
Our result can therefore be viewed as 
establishing an infinity-category of spectral sequences. And, via the shift-d\'ecalage adjunction, for each $r$ the relative category
$(\spse, \Ee_r)$ has $(\spse, \Ee_0)$ as a fibrant replacement.
\medskip

We note that there is some other recent work on spectral sequences, close in spirit to some of the ideas in this paper.
In addition to the previously mentioned work of Antieau~\cite{antieau2024}, in~\cite{blanc2022} the homotopy spectral sequence of
a (co)simplicial object in an infinity-category is studied and successive terms are analysed using two types
of localization. And~\cite{aitken24} studies
model structures on the infinity-category of filtered objects in a suitable stable infinity-category, noting that these
may be viewed as a categorification of the model structures for filtered complexes of~\cite{CELW20}. In~\cite{CLW25}, model category structures related to spectral sequences on the category of truncated multicomplexes are presented, together with potential applications to complex geometry.
\medskip

The paper is organized as follows. Section~\ref{sec:prelim} contains background material on presheaves, model categories
and bigraded complexes. In Section~\ref{sec:specseq}, we introduce the category of extended spectral sequences. The category $\D$ of disc objects is presented in Section~\ref{sec:lwb}, together with the presheaf category of linear witness books. 
Section~\ref{sec:qn} covers the main adjunction relating linear witness books and extended spectral sequences. In Section~\ref{S:TSD}
we study shift and d\'ecalage functors. The homotopical parts of the paper are in Section~\ref{S:ModelsESpSe}, on model category structures,
and in Section~\ref{S:infcat}, on an infinity-categorical interpretation. The paper concludes with an appendix of additional material which is helpful in providing a complete story of our presheaf approach, but not essential for the main text.
Appendix~\ref{App:liftingprops} contains details of important representable linear witness books. Appendix~\ref{App:adjoints} provides details of the left adjoint $\LDec$ to the shift functor.
And finally, Appendix~\ref{App:nonmodels} explains some non-existence results for certain model category structures having spectral sequences as fibrant objects.

\subsection*{Acknowledgements}
The authors would like to thank the Isaac Newton Institute for Mathematical Sciences, Cambridge, for support and hospitality during the programme \emph{Topology, representation theory and higher structures}, at Gaelic College, Sabhal Mòr Ostaig, Isle of Skye, where some of the work on this paper was undertaken. This work was supported by EPSRC grant EP/R014604/1.

\section{Preliminaries}
\label{sec:prelim}

Throughout this paper, we let $\kk$ denote a commutative ring with unit.  The category of $\kk$-modules is denoted $\modu$; it is a closed symmetric monoidal category.

\subsection{Linear presheaves}

We work with linear categories, that is, categories enriched in $\modu$, and we refer to~\cite{Kel05} for our conventions.
Given a small  linear category $\Bcat$ we denote by $\widehat{\Bcat}$ the category of  linear functors from 
${\Bcat}^{\mathrm{op}}$ to $\modu$. Such a category is called a category of {\sl linear presheaves}.
It is a  linear category with  linear natural transformations as morphisms. The Yoneda embedding $\mathcal 
Y_\Bcat\colon\Bcat\rightarrow \widehat{\Bcat}$ which associates to $b$ the functor 
$\Hom_{\Bcat}(-,b)$  is thus an enriched functor. We will denote the Yoneda functor $\mathcal 
Y_\Bcat$ simply by $\mathcal Y$ when the category $\Bcat$ is clear from the context.
We gather here the results in Section 4.1 of \cite{Kel05} that we will use in the paper.

\begin{prop}\label{P:Kelly} 
Let $\Bcat$ be a small  linear category and let $F\colon \Bcat \rightarrow \Ccat$ be an enriched functor. 
\begin{enumerate}
\item 
If $\Ccat$ is a cocomplete  linear category which is tensored 
over $\modu$, then the left Kan extension of $F$ along the Yoneda functor $
\mathrm{Lan}_{\Yy} F\colon \widehat{\Bcat}\rightarrow \Ccat$ exists and admits the following description as a coend. For $b'\in\widehat{\Bcat}$,
\[
	(\mathrm{Lan}_{\Yy} F)(b')=\int^{b\in \mathcal{B}} \Hom_{\widehat{\Bcat}}(\Yy(b),b')\otimes F(b).
\]
It is left adjoint to the functor  $\Ccat\rightarrow \widehat{\Bcat}$ which associates to $c$ the functor $\Hom_{\Ccat}(F(-),c)$.
\item If $\Ccat$ is a complete linear category which is cotensored 
over $\modu$, then the right Kan extension of $F$ along the Yoneda functor $
\mathrm{Ran}_{\Yy} F\colon \widehat{\Bcat}\rightarrow \Ccat$ exists and admits the following description as an end. 
 For $b'\in\widehat{\Bcat}$, 
\[
	(\mathrm{Ran}_{\Yy} F)(b')=\int_{b\in \mathcal{B}} \Hom\left( \Hom_{\widehat{\Bcat}}(b', \Yy(b)), F(b)\right).
\]
It is right adjoint to the functor $\Ccat\rightarrow \widehat{\Bcat}$ which associates to $c$ the functor $\Hom_{\Ccat}(F(-),c)$. \qed
\end{enumerate}
\end{prop}

The first adjunction corresponds to the usual adjunction between the 
nerve and the realisation functor.

Given a linear functor $\Bcat\to\Ccat$, applying the above to the composite with the Yoneda functor
$\mathcal 
Y_\Ccat\colon\Ccat\rightarrow \widehat{\Ccat}$ 
 yields the following.

\begin{prop}\label{P:presheaf-yoga}
A linear functor $G:\Bcat\to\Ccat$ gives rise to a linear functor $G^*\colon \widehat{\Bcat}\to \widehat{\Ccat}$ where
$G^*(b)=b\circ G$. This has a left adjoint $G_!$ and a right adjoint $G_*$. \qed
\end{prop}

\subsection{Model category structures}

We assume the reader is familiar with the language of model category structures and we will use standard model category theoretic
terminology and notation as can be found in~\cite{Hovey}. We will use standard recognition results for cofibrantly generated model
categories (cf.~\cite[Theorem 11.3.1]{Hir} or~\cite[Theorem 2.1.19]{Hovey}). 

We also use standard machinery for right transfer of model category structures, as summarised 
in~\cite[Section 4.4.1]{Balchin-handbook}. In particular, if we have an adjunction 
\[\xymatrix{
 \Ccat\ar@<1ex>[r]^-{L} & \Ccat' \ar@<1ex>[l]^-{R}_-{\perp} },
\]
with $\Ccat$ a model category, then we have a model category structure on $\Ccat'$, in which a map is a fibration or weak equivalence precisely if its image under $R$ is such in $\Ccat$, under either of the following conditions.
\begin{itemize}
    \item $\Ccat$ is cofibrantly generated, $R$ preserves filtered colimits, and maps in $\Ccat'$ with the left lifting property
    with respect to fibrations (that is, \emph{anodyne} maps) are weak equivalences. See~\cite[Proposition 3.1]{nlab:transferred_model_structure}.
    \item  $\Ccat$ is cofibrantly generated, $R$ also has a right adjoint $R'$
    and $(RL, RR')$ is a Quillen adjunction.
    See~\cite[Theorem 2.3]{DCH19}.
\end{itemize}

We make use of criteria for a Quillen adjunction to be a Quillen equivalence.
In particular, if in a Quillen adjunction 
 the right adjoint creates weak equivalences then the adjunction 
 is a Quillen equivalence precisely if for all cofibrant objects 
 the unit is a weak equivalence.
See~\cite[Proposition 2.3]{nlab:quillen_equivalence}.

We also use right Bousfield localization of a model category structure and the corresponding notion of right delocalization~\cite{C-S16}.

\subsection{Bigraded complexes}

In this section we let $r\geq 0$ be an integer.

\begin{defi}\label{D:bigraded modules} A \emph{bigraded $\kk$-module} $A$ is a collection of $\kk$-modules $A=\{A^{p,n}\}$ with $p,n\in\ZZ$. 
\end{defi}

\begin{defi}\label{D:r-complex} An \textit{$r$-bigraded complex}  is a bigraded $\kk$-module $A=\{A^{p,n}\}$ together with maps of $\kk$-modules
$d_r\colon A^{p,n}\to A^{p-r, n+1-r}$, called differentials, such that $d_r^2=0$. A \textit{morphism of $r$-bigraded complexes} is a map of bigraded modules commuting with the differentials. We denote by $\rbc$ the category of 
$r$-bigraded complexes. 
\end{defi}

\begin{rmk}
The grading convention here is chosen to be compatible with our previous work~\cite{CELW20,LW24}. It will mean that the differential on the $r$-page
of a spectral sequence has bidegree $(-r, 1-r)$. It is straightforward to translate our results to the standard setting of a homological spectral sequence where the corresponding bidegree is $(-r, r-1)$ or that of a cohomological spectral sequence where it is $(r,1-r)$.
\end{rmk}

\begin{notation}
We denote by $B_r(A)$ the image of $d_r$, that is the bigraded $\kk$-module of $r$-boundaries of $A$, and by $Z_r(A)$ the kernel of $d_r$, that is, the  bigraded $\kk$-module of $r$-cycles of $A$.  We write $H(A)$ for the homology $Z_r(A)/B_r(A)$.  
\end{notation}

We note that $\rbc$ can be viewed as a category of linear presheaves as follows.

\begin{defi}\label{D: delrcat}
For $r\geq 0$, we let $\delrcat$ be the following linear category. The objects are $(p,n)$ for $p,n\in\ZZ$. The $\kk$-module of morphisms $\Hom_{\delrcat}((p,n), (q,m))$ is either $0$ or free of rank $1$ spanned by $\delta_r^{p,n}\colon (p-r,n-r+1)\rightarrow (p,n)$ or by the identity. Hence, composition of two non-identity morphisms is $0$. 
\end{defi}

There is an isomorphism of categories 
\begin{equation}\label{iso:rC} \widehat{\delrcat}\cong \rbc.\end{equation}
 On objects this is given by the assignment that sends $A$ in $\widehat{\delrcat}$ to the $r$-bigraded complex with $A^{p,n}=A(p,n)$ and with $d_r\colon A^{p,n}\to
 A^{p-r,n+1-r}$ given by $A(\delta_r^{p,n})$.

\begin{notation}
We write $\kk^{p,n}$ for the bigraded module consisting of a free module of rank $1$ in bidegree $(p,n)$ and
otherwise $0$.
\end{notation}

\begin{rmk}
Note that the Yoneda functor 
$\Yy\colon \delrcat\rightarrow \widehat{\delrcat}$ sends $(p,n)$ to the $r$-bigraded complex  having only two non-zero components, each free of rank 1, at bigrading $(p,n)$ and $(p-r,n+1-r)$ with $d_r\colon \kk\rightarrow\kk$ the identity map.
We have $\kk^{p-r,n-r+1}=\mathrm{coker}( \delta_r\colon
\Yy(p-2r,n+2-2r)\rightarrow \Yy(p-r,n+1-r) )$ and since $\delta_r^2=0$ we have
an induced map $\delta_r^{p,n}\colon \kk^{p-r,n+1-r}\rightarrow \Yy(p,n)$.
\end{rmk}

The homology of an $r$-bigraded complex  is a bigraded $\kk$-module and the category of $r$-bigraded complexes has a natural class of quasi-isomorphisms, namely morphisms inducing isomorphisms on homology.

Similarly to the treatment of unbounded chain complexes by Hovey in ~\cite[Section 2.3]{Hovey}, the $r$-bigraded complex $\Yy(p,q)$ corresponds to the disc object $D^*(\kk)$, the $r$-bigraded complex $\kk^{p,n}$ to the sphere object $S^*(\kk)$, and the map $\delta^r: \kk^{p-r,n+1-r}\to\Yy(p,n)$ to the injection $S^{*-1}(\kk)\to D^*(\kk)$. 
The proposition below provides
 a model category structure corresponding to the projective model category structure on the category of unbounded chain complexes adapted to the category of $r$-bigraded complexes.

\begin{prop}\label{P:modcatrbc} The category $\rbc$ of $r$-bigraded complexes has a cofibrantly generated model category structure where
\begin{itemize}
\item fibrations are bidegreewise surjections,
\item weak equivalences are quasi-isomorphisms.
\end{itemize}
 The set of generating acyclic cofibrations is the set  $I=\{0\rightarrow \Yy(p,n)\}_{p,n\in\ZZ}$ and the set of generating cofibrations is the set 
  $J=\{\delta_r^{p,n}\colon \kk^{p-r,n-r+1}\rightarrow \Yy(p,n)\}_{p,n\in\ZZ}$. \qed
\end{prop}

We will identify $\widehat{\delrcat}$ with $\rbc$ via the isomorphism (\ref{iso:rC}) and view it as a model category with the above structure, which we refer to as the projective model structure.
\smallskip

We end this section with the $r$-cone construction that will be useful in the sequel.

\begin{defi}\label{D:rcone1}
Let $A$ be an $r$-bigraded complex. The \textit{$r$-cone of $A$}
is the $r$-bigraded complex $(\cone_r(A),d_r)$ given by
\[\cone_r(A)^{p,n}=A^{p,n} \oplus A^{p+r,n+r-1}\text{ with }d_r(a,b)=(0,a).\]
It is naturally endowed with a projection $\pi_r\colon \cone_r(A)\rightarrow A$ defined as $\pi_r(a,b)=a+d_r^A(b)$.
\end{defi}

Note that the assignment $A\mapsto \cone_r(A)$ is functorial, and that this functor is the composite of the forgetful functor from  $\rbc$ to $\bimod$  followed by its left adjoint. Note also that it is isomorphic to the usual cone construction of the identity map in chain complexes. In particular, for every $A$, $\cone_r(A)$ is acyclic. 

\smallskip

We end this subsection by noting that we can compare 
 the categories of $r$-bigraded complexes for varying $r$ by using a translation functor.

\begin{defi}\label{Def:shift_rbc}
For each $r\geq 0$, the \emph{translation functor} $\Tt\colon  \rbc \to (r+1)\text{-}\mathsf{C}$ is given by
 $(\Tt A)^{p,n}=A^{n, 2n-p}$, with differential that of $A$, and $(\Tt f)^{p,n}=f^{n, 2n-p}$.
\end{defi}

Note that translation commutes with cone constructions and homology.

\begin{prop}\label{P:Qequiv_cxs}
For each $r\geq 0$, 
the categories of $r$-bigraded complexes and $(r+1)$-bigraded complexes are isomorphic via the translation functor,
$\Tt\colon \rbc \to (r+1)\text{-}\mathsf{C}$. This induces a Quillen equivalence between the projective model category structures on $\rbc$ and $(r+1)\text{-}\mathsf{C}$. \qed
\end{prop}

\section{The category of extended spectral sequences}\label{sec:specseq}

In~\cite{LW24} we studied the category of spectral sequences $\spse$, and provided some examples showing that it is neither complete nor cocomplete. In this section we embed the category $\spse$ in the category of extended spectral sequences. We show here that this category is cocomplete and we will see later that it is also complete; see Proposition~\ref{P:lwbe-lims-colims}.

\subsection{Spectral sequences and extended spectral sequences}\label{subsec:defs}

\begin{defi}\label{def:spse_obj}
An \textit{extended spectral sequence} $(X,\varphi)$ is a family of $r$-bigraded complexes $(X_r, d_r)$, for $r\geq 0$, together with a family of morphisms of bigraded $\kk$-modules $\varphi_{r+1}\colon  X_{r+1} \rightarrow H(X_r)$ for $r\geq 0$, called \emph{characteristic maps}.

A \textit{morphism of  extended spectral sequences} is a family of morphisms $f_r\colon X_r\to Y_r$ of $r$-bigraded complexes, for $r\geq 0$,
which is \textit{compatible with characteristic maps}. 

We denote by $\espse$ the category of extended spectral sequences. 

A \emph{spectral sequence} is an extended spectral sequence such that the characteristic maps are all isomorphisms. We denote by $\spse$ the category of spectral sequences, that is, the full 
subcategory of $\espse$ whose objects are spectral sequences. 
\end{defi}

We note that the categories $\spse$ and $\espse$ implicitly involve the
forgetful functors from $r$-bigraded complexes to bigraded $\kk$-modules,
since the characteristic maps are morphisms of  
bigraded $\kk$-modules.

\begin{prop}\label{P:espse_colimit} The category $\espse$ of extended spectral sequences is a linear category with (small) colimits,  computed pagewise. It is tensored over $\modu$.
\end{prop}

\begin{proof} That the category is linear is immediate. It is straightforward to see that it is tensored over $\modu$, noting that, for 
an $\kk$-module $M$ and an extended spectral sequence $Y$, the characteristic maps of $M\otimes Y$ are given by the composites
\[
	M\otimes Y_{i+1}\xrightarrow{1_M\otimes \varphi^Y_{i+1}} M\otimes H(Y_i)\xrightarrow{\mu}  H(M\otimes Y_i),
\]
where $\mu(m\otimes[x])=[m\otimes x]$.

 Let $X\colon I\rightarrow \espse$ be a functor with $I$ a (small) category. 
For each $m\geq 0$ denote by 
\[
	\varphi(i)_{m+1}\colon X(i)_{m+1}\rightarrow H(X(i)_m)
\]
the characteristic maps, and by 
$Y_m=\mathrm{colim}_{i\in I} X(i)_m$ the colimit of the  diagram $X_m\colon I\rightarrow \mbc$ induced by the projection of $\espse$ to $\mbc$. The canonical maps $\rho(i)_{m}\colon  X(i)_{m}\rightarrow Y_m$ induce maps of bigraded $R$-modules
\[
	H(\rho(i)_{m})\circ \varphi(i)_{m+1}\colon X(i)_{m+1}\rightarrow H(Y_m).
\]
Since the forgetful functor from  $(m+1)$-bigraded complexes to bigraded $R$-modules preserves colimits, there is
 a (unique) characteristic map 
$Y_{m+1}\rightarrow H(Y_m)$. Thus the pagewise colimit $Y$ is an object of $\espse$, with compatible maps $X(i)\to Y$. For the universal property, if $Y'$ is another such with canonical maps $\rho'(i)_m\colon X(i)_m\to Y'_m$, there is a unique map of $m$-bigraded complexes $\alpha_m\colon Y_m\to Y'_m$ for all $m$, such that $\alpha_m\rho(i)_m=\rho'(i)_m$. And these satisfy compatibility with the characteristic maps by using the universal property of (the underlying bigraded module of) the pagewise colimit.
\end{proof}

\begin{rmk}
Note that $\espse$ is not cotensored over $\modu$. This is one consequence of the choice of direction of the characteristic maps.
\end{rmk}

\subsection{Discs in extended spectral sequences}\label{SS:disc}
We adapt the notation of Section 5.5 in~\cite{LW24} to the case of extended spectral sequences.

\begin{defi}\label{def:discs} 
Let $p,n\in\ZZ$. For all $r\geq 0$, let $\Dd_r(p,n)$ be the spectral sequence defined as follows.
\[
\begin{cases}
\Dd_r(p,n)_i=\kk^{p,n}\oplus \kk^{p-r,n+1-r}, \quad d_i=0 &\text{for } 0\leq i<r, \\
\Dd_r(p,n)_r=\kk^{p,n}\stackrel{1}{\lra}\kk^{p-r,n+1-r}\\
\Dd_r(p,n)_i=0 &\text{for } i>r.
\end{cases}
\]
We denote by $e_{r,i}^{p,n}$ a generator of $\Dd_r(p,n)_i^{p,n}$ and by $f_{r,i}^{p-r,n+1-r}$ a generator of 
$\Dd_r(p,n)_i^{p-r,n+1-r}$, for $0\leq i\leq r$. The characteristic map 
\[
    \varphi_i\colon  \Dd_r(p,n)_i\to H( \Dd_r(p,n)_{i-1}) 
\]
is given by
\begin{align*}
\varphi_i(e_{r,i}^{p,n})&=[e_{r,i-1}^{p,n}],\\
\varphi_i(f_{r,i}^{p-r,n+1-r})&=[f_{r,i-1}^{p-r,n+1-r}].
\end{align*}
for $0\leq i\leq r$ and $\varphi_i=0$ for $i>r$.

\end{defi}

\begin{defi}\label{D:representables} Let $(X,\varphi)$ be an extended spectral sequence. 
\begin{enumerate}
\item A sequence of elements $(x_0^{p,n},\ldots,x_{m+1}^{p,n})$ with $x_i^{p,n}\in X_i^{p,n}$ is said to be 
\emph{compatible} if for every $0\leq i\leq m$, $d_ix_i^{p,n}=0$ and $\varphi(x_{i+1}^{p,n})=[x_i^{p,n}]$ where $[x_i^{p,n}]$ is the class of $x_i^{p,n}$ in $H(X_i)$.
\item Denote by $\Nn(X)_r^{p,n}$ the $\kk$-submodule of $X^{\oplus (2r+2)}$ consisting of pairs
\[(x_0^{p,n},\ldots,x_r^{p,n});(y_0^{p-r,n+1-r},\ldots,y_r^{p-r,n+1-r})\] of compatible sequences satisfying $d_rx_r^{p,n}=y_r^{p-r,n+1-r}$. This yields a functor, denoted by $\Nn_r^{p,n}\colon \espse\rightarrow \modu$.
\end{enumerate}
\end{defi}

\begin{notation}
We denote by $\kk(p,n)$ the spectral sequence which at each page has a free module of rank $1$ concentrated in bidegree $(p,n)$ and all differentials zero.
\end{notation}

The following proposition is a direct consequence of the definitions.

\begin{prop}\label{P:morphismfromdisks} Let $(X,\varphi)$ be an extended spectral sequence.
\begin{enumerate}
\item There is a one-to-one correspondence between infinite compatible sequences $(x_0^{p,n},\ldots,x_m^{p,n},\ldots)$ and morphisms of extended spectral sequences $\kk(p,n)\rightarrow X$.
\item We have $\Nn_r^{p,n}=\Hom_{\espse}(\Dd_r(p,n),-)$, that is, $\Nn_r^{p,n}$ is represented by $\Dd_r(p,n)$.\qed
\end{enumerate} 
\end{prop}

\begin{rmk}
Later we will see that the functors $\Nn_r^{p,n}$ assemble to a nerve functor from extended spectral sequences to a presheaf category; see
Proposition~\ref{P:nerve_descr}.
\end{rmk}

\section{The presheaf category of linear witness books}
\label{sec:lwb}

In this section we introduce and study a linear presheaf category closely related to spectral sequences.

\subsection{The category \(\D\)}
We define the underlying category on which we will consider presheaves. We show in Proposition~\ref{prop:subcatondiscs}
that it can 
 be understood as the full linear subcategory of $\spse$ generated by the disc objects $\Dd_r(p,n)$. 

\begin{defi}\label{def:S}
Let $\D$ be the small linear category defined as follows. Objects of $\D$ are triples $(r,p,n)$ with $r\geq 0$ and $p, n\in \ZZ$.
Morphisms in the category are generated by the following  morphisms
for all $r\geq 0$ and $p,n\in\ZZ$.
\begin{align*}
\textit{co-witness}\qquad &\omega_{r+1}^{p,n}\colon (r,p,n)\rightarrow (r+1,p,n), \\
\textit{co-differential}\qquad &\delta_r^{p,n}\colon (r,p-r,n-r+1)\rightarrow (r,p,n),\\
\textit{co-degeneracy}\qquad &\sigma_r^{p,n}\colon (r+1,p+1,n+1)\rightarrow (r,p,n), 
\end{align*}
subject to the following relations for all $r\geq 0$ and  $p,n\in\ZZ$.
\[
\sigma_{r}^{p-1,n-1}\omega_{r+1}^{p,n}=
0, \qquad
\omega_{r+1}^{p,n}\sigma_{r}^{p,n}=0,  \qquad
\sigma_{r}^{p+r,n+r-1} \delta_{r+1}^{p+r+1,n+r} \omega_{r+1}^{p,n}=\delta_r^{p+r,n+r-1}.
\]
\end{defi}

\begin{rmk}
The relations
\[
\delta_r^{p+r,n+r-1}\delta_r^{p,n}=\delta_{r}^{p+r,n+r-1}\sigma_r^{p,n}=\omega_{r+1}^{p,n}\delta_r^{p,n}=0
 \qquad \text{ for all } r\geq 0, p,n\in\ZZ
\]
also hold in $\D$ as an immediate consequence.
\end{rmk}

The following lemma is a direct consequence of the definition and gives all the non-zero $\kk$-modules of morphisms in $\D$.
\begin{lem}\label{L:morinS}
For every $p,n\in\ZZ$, $r\geq 0$,
\begin{align*}
&\Hom_{\D}((r-i,p,n),(r,p,n))&&= \kk (\omega_r)^i&\ \text{for }  r\geq i\geq 0,\\
&\Hom_{\D}((r+j,p+j,n+j),(r,p,n))&&= \kk (\sigma_{r})^j&\ \text{for }  j>0,\\
&\Hom_{\D}((r,p-r,n-r+1),(r,p,n))&&= \kk \delta_r, &\\
&\Hom_{\D}((r-i,p,n),(r,p+r,n+r-1))&&= \kk \delta_{r}(\omega_r)^i&\ \text{for }  r\geq  i>0,\\
&\Hom_{\D}((r+j,p,n),(r,p+r,n+r-1))&&= \kk(\sigma_{r})^j \delta_{r+j}&\ \text{for }  j>0.
\end{align*}
Here we have omitted the upper indices on the morphisms; they can be determined by the hom modules in which the morphisms live. 
We denote by $(\omega_r)^i$ the composite $\omega_{r}\omega_{r-1}\ldots \omega_{r-i+1}$ and the same notation is used for the 
iterated composite of the $\sigma$s. For $i=0$, it denotes the identity map.
\qed
\end{lem}

\begin{rmk}
     The category $\D$ is self-dual. Indeed,  
     for a morphism $\alpha\colon C\to C'$ in a category $\Ccat$ we
write $\overline{\alpha}$ for the corresponding morphism $C'\to C$ in $\Ccat^{\mathrm{op}}$.
Then
     an isomorphism $\D \to \D^{\mathrm{op}}$ is given on objects
    by $(r,p,n)\mapsto (r,r-p,r-n)$ and on generating
    morphisms by $
        \delta_r^{p,n}\mapsto \overline{\delta_r^{2r-p, 2r-n-1}},
        \sigma_r^{p,n} \mapsto \overline{\omega_{r+1}^{r-p,r-n}},
        \omega_{r+1}^{p,n}\mapsto \overline{\sigma_r^{r-p, r-n}}.
   $
\end{rmk}

We define some maps of spectral sequences between the disc objects $\Dd_r(p,n)$ of Definition~\ref{def:discs}.

\begin{defi}\label{def:mapsofdiscs}
Recall that we denote by $e_{r,i}^{p,n}$ a generator of $\Dd_r(p,n)_i^{p,n}$ and by $f_{r,i}^{p-r,n+1-r}$ a generator of 
$\Dd_r(p,n)_i^{p-r,n+1-r}$, for $0\leq i\leq r$. We define maps of spectral sequences 
\begin{align*}
&\omega\colon \Dd_r(p,n)\to \Dd_{r+1}(p,n)\\
&\delta\colon \Dd_r(p,n)\to \Dd_{r}(p+r,n+r-1)\\
&\sigma\colon \Dd_r(p,n)\to \Dd_{r-1}(p-1,n-1)
\end{align*}
determined by, for $0\leq i\leq r$, 
\begin{align*} 
\omega(e_{r,i}^{p,n})&=e_{r+1,i}^{p,n},  \\
 \omega(f_{r,i}^{p-r,n+1-r})&=0,\\
\delta(e_{r,i}^{p,n})&=f_{r,i}^{p,n},\\
\delta(f_{r,i}^{p-r,n+1-r})&=0,\\
 \sigma(e_{r,i}^{p+1,n+1})&=0,
 \end{align*}
 and
\begin{align*} 
\sigma(f_{r,i}^{p-r,n+1-r})&=f_{r-1,i}^{p-r,n+1-r}\quad \text{for }0\leq i< r,\\
\sigma(f_{r,r}^{p-r,n+1-r})&=0.
\end{align*}
\end{defi}

\begin{rmk}
If we think of the generators $e$ and $f$ as the \emph{source} and \emph{target} respectively, $\omega$ maps source to source, $\delta$ maps source to target and $\sigma$ maps target to target. 
\end{rmk}

In the following diagrams both  $\kkpic$ and $\kkpnpic$ represent $\kk$, a free module of rank $1$, the vertex marked $\kkpnpic$ being in bidegree $(p,n)$. The solid arrows are identity maps within disc objects and the dashed arrows are identity maps between source or target of the disc objects as indicated.

The morphism $\omega\colon \Dd_r(p,n)\to \Dd_{r+1}(p,n)$ can be depicted as follows.

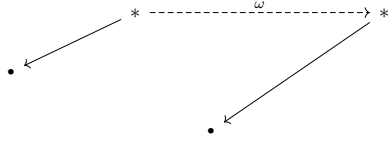
\vspace{0.8cm}
  \begin{center}  
    \begin{tikzcd}[transform canvas={scale=0.7}]
      && \kkpnpic   \arrow[lld]   \arrow[rrrr, dashed]{}{\omega}&&&&\kkpnpic \arrow[llldd]\\
      \kkpic&&&&&&\\
      &&&\kkpic &&&                             
    \end{tikzcd}
  \end{center}
  \vspace{0.8cm}

The corresponding picture for $\delta\colon \Dd_r(p,n)\to \Dd_{r}(p+r,n+r-1)$ is:
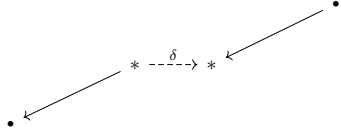
\vspace{0.8cm}
  \begin{center}  
    \begin{tikzcd}[transform canvas={scale=0.7}]
      &&&&& \kkpic \arrow[lld]\\
      &&\kkpnpic \arrow[lld]\arrow[r,dashed]{}{\delta}&\kkpnpic&&\\
      \kkpic &&&&&                             
    \end{tikzcd}
  \end{center}
  \vspace{0.8cm}

and for $\sigma\colon \Dd_r(p,n)\to \Dd_{r-1}(p-1,n-1)$:
\vspace{0.8cm}
  \begin{center}  
    \begin{tikzcd}[transform canvas={scale=0.7}]
      &&& \kkpnpic \arrow[llldd]&&&\\
      &&&&&&\kkpic \arrow[lld]\\
      \kkpic\arrow[rrrr,dashed]{}{\sigma} &&&&\kkpic&&                             
    \end{tikzcd}
  \end{center}
  \vspace{0.8cm}

\begin{prop}
\label{prop:subcatondiscs}
 The category $\D$ is isomorphic to the full subcategory of $\spse$ generated by the objects $\Dd_r(p,n)$ for $r\geq 0$ and $p,n\in\ZZ$, via the assignment $(r,p,n)\mapsto \Dd_r(p,n)$.
\end{prop}

\begin{proof} 
We send the generating morphisms of $\D$ to the correspondingly named morphisms of spectral sequences given in Definition~\ref{def:mapsofdiscs}. It is straightforward to check that the relations of Definition~\ref{def:S} are satisfied by these morphisms. Thus we have a well-defined functor from $\D$ to $\spse$.

We claim that non-trivial morphisms of spectral sequences $\Dd_r(p,n)\rightarrow \Dd_{r'}(p',n')$ correspond
exactly to the non-trivial morphisms in $\D$, as specified in Lemma~\ref{L:morinS}.

This is checked by direct inspection of morphisms of spectral sequences $\Dd_r(p,n)\rightarrow \Dd_{r'}(p',n')$. Briefly, we can have such a morphism of type ``source to source'' only for $r<r'$ and these are generated by composites of $\omega$s. There is a morphism of type ``target to target'' only for $r>r'$ and these are generated by composites of $\sigma$s. And morphisms of type
``source to target'' are generated by $\delta$ if $r=r'$, by composites of the form $\delta\omega^i$ if $r<r'$ and
by composites of the form $\sigma^j\delta$ if $r>r'$. There is no non-trivial map of spectral sequences of type ``target to source'' and, together with identity maps, this accounts for all the morphisms.
\end{proof}

\subsection{Linear witness books}

\begin{defi} We write $\lwb$ for the linear presheaf category  $\widehat{\D}$. 
 A \emph{linear witness book} is an object of $\widehat{\D}$ and
a \emph{morphism of  linear witness books} is a morphism in   $\widehat{\D}$.
\end{defi}

The following proposition is a direct consequence of the definitions.
\begin{prop}\label{P:hats}
A linear witness book $L$ is a collection of $r$-bigraded complexes $(L_r,d_r)$ for $r\geq 0$, endowed with
$\kk$-linear witness maps $w_{r+1}\colon L_{r+1}\rightarrow L_{r}$ of bidegree $(0,0)$ and $\kk$-linear degeneracy maps $s_r\colon L_r\rightarrow L_{r+1}$ of bidegree $(1,1)$ satisfying the relations
	\[
	w_{r+1}s_{r}=0,\; s_rw_{r+1}=0 \text{ and  }d_r=w_{r+1}d_{r+1}s_r,\ \text{ for all } r\geq 0.
	\]
A morphism of linear witness books $f\colon K\to L$ is a collection of morphisms of $r$-bigraded complexes
$f_r\colon K_r\to L_r$ for $r\geq 0$, compatible with the
witness and degeneracy maps.
\qed
\end{prop}

\begin{defi}
The \emph{$r$-page} of a linear witness book $L\in\lwb$ is the  $r$-bigraded complex $(L_r,d_r)$.
\end{defi}

We write $(L,d_i,w_i,s_i)$ for an object of $\lwb$.

\begin{rmk}
For $(L,d_i,w_i,s_i)\in\lwb$, we also have the relations $s_rd_r=d_rw_{r+1}=0$. 
The relation
$d_r^2=0$ is redundant, but the connection to spectral sequences is clearer if we explicitly view the $r$-page as an $r$-bigraded complex.
\end{rmk}

It follows directly from the description above that for any linear witness book $(L, d_i,w_i,s_i)$ we have, for all $i\geq 0$,
\begin{equation}\label{E:fundamental}
\Img d_i\subseteq \Img w_{i+1}\subseteq \Ker s_i\subseteq \Ker d_i \ \text{ and }\  \Img s_i\subseteq \Ker w_{i+1}.
\end{equation}

The combinatorics of the representable object $\Yy(r,p,n)$ are treated in Example~\ref{E:Yrpn}.

\begin{prop}\label{P:QN} The embedding $\D\hookrightarrow \espse$ which associates to $(r,p,n)$ the spectral sequence $\Dd_r(p,n)$ induces an adjunction
\[\xymatrix{
 \lwb\ar@<1ex>[rr]^-{\Qq} && \espse \ar@<1ex>[ll]^-{\Nn}_-{\perp} }\]
 with the right adjoint defined as $\Nn(X)(r,p,n)=\Hom_{\espse}(\Dd_r(p,n),X).$
\end{prop}

\begin{proof} 
Apply Proposition~\ref{P:Kelly}, noting that $\espse$ is a cocomplete linear category tensored over $\modu$ by Proposition~\ref{P:espse_colimit}.
\end{proof}

By Proposition~\ref{P:morphismfromdisks}, the functor $\Nn$ is explicitly described in
Definition~\ref{D:representables}.

\section{Relating linear witness books and extended spectral sequences}
\label{sec:qn}

In this section we explore how the categories of linear witness books and extended spectral sequences are related. We begin by supplying more details of the 
$(\Qq, \Nn)$ adjunction and explaining its properties. This will allow us to 
identify a full subcategory of $\lwb$ that is equivalent to $\espse$. Importantly,
this is a reflective subcategory, allowing us to conclude in 
Theorem~\ref{T:espse-lims-colims} that $\espse$ is bicomplete. We also provide characterizations of (extended) spectral sequences inside $\lwb$.

\subsection{Description and properties of the main adjunction}

Though the $(\Qq, \Nn)$ adjunction comes from standard category theory, it is important to describe the functors. The functor $\Nn$, similar to a nerve functor, is easily defined. The functor $\Qq$, obtained via a coend formula, has a nice description as a quotient by degeneracies.

\begin{prop}\label{P:nerve_descr} For the functor $\Nn\colon \espse\rightarrow \lwb$, the module
$\Nn(X)_r^{p,n}$ is as described explicitly in part (2) of Definition~\ref{D:representables}.
In addition,
 for any pair of compatible sequences $(x;x')=((x_0,\ldots,x_r);(x'_0,\ldots,x'_r))$  in $\Nn(X)_r$ with $d^X_rx_r=x'_r$ we have
\[d_r(x;x')=(x';0),\; w_r(x;x')=((x_0,\ldots,x_{r-1});0),\; s_r(x;x')=(0;(x'_0,\ldots,x'_r,0)).\] \end{prop}
\begin{proof}
We have to show that the maps $d_r,w_r,s_r$ have the correct expression, which is a matter of directly interpreting the expression $\Nn(X)_r^{p,n}=\Hom_{\espse}(\Dd_r(p,n),X).$
\end{proof}

\begin{notation}\label{N:S_r}
For a linear witness book $L$ and $r>0$, we write
$S_r=S_r(L)$ for the sub $r$-bigraded complex $s_{r-1}L_{r-1}+d_rs_{r-1}L_{r-1}$ of $L_r$.
\end{notation}

\begin{teo}\label{T:adj_QN}

The functor $\Qq\colon \lwb\rightarrow \espse$ assigns to $L\in\lwb$ the extended spectral sequence $\Qq(L)$, where $\Qq(L)_0=L_0$ and for $r>0$,
\[\Qq(L)_r=L_r/S_r,\]
with differential induced by that of $L_r$
and with characteristic map $\varphi_r\colon \Qq(L)_r\rightarrow H(\Qq(L)_{r-1})$ defined by
\[\varphi_r(\overline a)=\left[\overline{w_r(a)}\right].\]
Here bars denote classes in the quotient by $S_r$ and square brackets denote homology classes. \\
For $X$ an extended spectral sequence, the counit of the $(\Qq,\Nn)$  adjunction
 \[\epsilon_X\colon \Qq\Nn(X)\rightarrow X\] is induced by
the projection $\Nn(X)_r\rightarrow X_r$ which associates $x_r$ to $(x;x')$. \\
The counit is a componentwise isomorphism.
 In particular, $\Nn$ is fully faithful and $\Qq$ is essentially surjective.\\
For $L$ a linear witness book,
the unit of the  $(\Qq,\Nn)$ adjunction 
\[\eta_L\colon L\rightarrow \Nn\Qq(L),\]
 is given, for $a\in L_r$, by 
	\[
	(\eta_L)_r(a)=\left(\overline{w^r(a)},\ldots,\overline{w(a)},\overline{a};
                \overline{w^r(d_ra)},\ldots,\overline{w(d_ra)},\overline{d_ra}\right).
	\]
	  The unit is componentwise surjective.
\end{teo}

\begin{proof} Though the description of $\Qq$ has been obtained by the coend formula of Proposition \ref{P:Kelly}, we prove that $(\Qq, \Nn)$ is an adjunction by checking the triangle identities relating the unit and the counit. That $\epsilon$ and $\eta$ are well defined natural transformations is an easy check. We fix $X$ an extended spectral sequence and $L$ a linear witness book.
Proposition \ref{P:nerve_descr} implies that 
\[S_r(\Nn(X))=\{(x;x')\in\Nn_r(X)\,|\, x_r=0\}.
\] In particular, in $\Qq\Nn(X)_r$, one has $\overline{(x;x')}=\overline{(y;y')}$ if and only if $x_r=y_r$. We claim that this implies the following.
\begin{itemize}
    \item[a)] $(\Qq\eta_L)\circ \epsilon_{\Qq L}=1_{\Qq\Nn\Qq L}.$
    \item[b)]  The counit $\epsilon_X$ is componentwise injective.
     \item[c)]  $(\eta_{\Nn X})\circ \Nn\epsilon_{X}=1_{\Nn\Qq\Nn X}.$
\end{itemize}
For a) and b) this is a direct check.
Let us prove c). Let $(x;y)\in (\Nn\Qq\Nn X)_r$ with $x=(\bar x_0,\ldots,\bar x_r)$ and $y=(\bar y_0,\ldots,\bar y_r)$ where $x_i$ and $y_i$ are in $(\Nn X)_i$. Denote by $x_i^i$ and $y_i^i$ their projections onto $X_i$, so that $a:=\Nn\epsilon_X(x;y)=(x_0^0,\ldots,x_r^r;y_0^0,\ldots,y_r^r)\in (\Nn X)_r$.
We have, for $i\geq 1$,
$w^i(a)=(x_0^0,\ldots,x_{r-i}^{r-i};0)$ and 
$w^i d_ra=(y_0^0,\ldots,y_{r-i}^{r-i},0)\in(\Nn X)_{r-i}$. As a consequence $\overline{w^i a}=\bar x_{r-i}$ and $\overline{w^i d_ra}=\bar y_{r-i}$ in $(\Qq\Nn X)_{r-i}$. This shows that $\eta_{\Nn X}(a)=(x;y)$.
\\
Claims a) and c) imply that the left adjoint of $\Nn$ is the described functor $\Qq$. \\
The counit is componentwise surjective because the projection $\Nn(X)_r\to X_r$ is. 
Let us prove that the unit is componentwise surjective.\\
Let $\alpha=(\overline{a_0}, \dots, \overline{a_r}; \overline{b_0}, \dots, \overline{b_r})\in
 \Nn\Qq(L)_r$. We claim that for $j=0, \dots , r$ there is an element $c_{r-j}\in L_r$ such that 
 $\overline{w^i c_{r-j}}=\overline{a_{r-i}}$ and $\overline{w^i d_r c_{r-j}}=\overline{b_{r-i}}$ for $0\leq i\leq j$.  Then $(\eta_L)_r(c_0)=\alpha$ and $(\eta_L)_r$
 is surjective.
 
 We prove the claim by induction on $j$. For $j=0$, we can take $c_r=a_r$. Now assume we have $c_{r-j}$ as above. We have
 $\left[\overline{w^{j+1}c_{r-j}}\right]=[\overline{a_{r-j-1}}]$ in $H((\Qq L)_{r-j-1})$. So there is some $c'\in L_{r-j-1}$ such that
 \[\overline{w^{j+1}c_{r-j}}=\overline{a_{r-j-1}} +\overline{d_{r-j-1}}\overline{c'}
 =\overline{a_{r-j-1}} +\overline{w_{r-j}d_{r-j}s_{r-j-1}c'}\]
 and so
 \[
 \overline{w_{r-j}(w^jc_{r-j}-d_{r-j}s_{r-j-1}c')}=\overline{a_{r-j-1}}.
 \]
 Similarly, since  $[\overline{w^{j+1}d_rc_{r-j}}]=[\overline{b_{r-j-1}}]$ in $H((\Qq L)_{r-j-1})$, there is some $c''\in X_{r-j-1}$ such that
  \[
 \overline{w_{r-j}(w^jd_rc_{r-j}-d_{r-j}s_{r-j-1}c'')}=\overline{b_{r-j-1}}.
 \]
 Then it is straightforward to check that $c_{r-j-1}=c_{r-j}-d_rs^{j+1}c'-s^{j+1}c''$ has the required properties.
\end{proof}

\begin{prop}\label{P:projection}
There is a forgetful functor $\Uu\colon \lwb\to\espse$ and the quotient maps
$L_r\to \Qq(L)_r$ give rise to a natural transformation of functors from $\Uu$ to $\Qq$.
\end{prop}

\begin{proof}
Let $L$ be a linear witness book.
Since $d_{r-1}w_r=0$,
the collection $(L_r, d_r)$ of $r$-bigraded complexes may be endowed with the characteristic maps
$\varphi_r\colon L_r\rightarrow H(L_{r-1})$ given by
$\varphi_r(a)=[w_r(a)]$. This gives the extended spectral sequence $\Uu(L)$. It is clear that a morphism of
linear witness books $f\colon L \to K$ gives a morphism of extended spectral sequences $\Uu(L)\to \Uu(K)$.

By Theorem~\ref{T:adj_QN} $\Qq(L)_r$  is isomorphic to the quotient 
$L_r/S_r$. The description of 
 the characteristic maps of $\Qq(L)$ in Theorem~\ref{T:adj_QN} shows that this gives a 
 morphism  $\Uu(L)\to\Qq(L)$ in $\espse$.
\end{proof}

The following result will be important for understanding model category structures later.

\begin{prop}\label{P:Nsurj}
Let $f\colon X\to Y$ be a morphism in $\espse$. Then $f_i$ is surjective for $0\leq i\leq r$ if and only if
$\Nn(f)_i$ is surjective for $0\leq i\leq r$.
\end{prop}

\begin{proof} It is clear that if $\Nn(f)_r$ is surjective then $f_r$ is surjective.
Let us prove that if $f_i$ is surjective for $0\leq i\leq r$, then so is $\Nn(f)_i$ for $0\leq i\leq r$, by induction on $r$. This is true for $r=0$.
Let $X$ be an extended spectral sequence.
We have seen in the proof of  Theorem~\ref{T:adj_QN} that Proposition~\ref{P:nerve_descr} implies that $S_r(\Nn X)$ is the kernel of the projection $(\Nn X)_r\to X_r$. In addition,
$S_r(\Nn X)=s_{r-1}( \Nn(X)_{r-1})\oplus  d_rs_{r-1}(\Nn(X)_{r-1}).$
Hence, for 
$r\geq 1$ and $p,n\in\ZZ$, there is a natural short exact sequence of $\kk$-modules
\[\xymatrix{
s_{r-1}( \Nn(X)_{r-1}^{p-1,n-1})\oplus  d_rs_{r-1}(\Nn(X)_{r-1}^{p+r-1,n+r-2})\ar@{^{(}->}[r]^-{}
& \Nn(X)_r^{p,n}\ar@{-{>>}}[r]& X_r^{p,n}
 }\]
This implies that if $\Nn(f)_{r-1}$ is surjective, then $f_r$ is surjective if and only if $\Nn(f)_r$ is surjective. 
\end{proof}

\subsection{Characterizing extended spectral sequences as presheaves}

We characterize the essential image of the functor $\Nn$, thus giving a subcategory of the category of linear
witness books $\lwb$ which is equivalent to the category of extended spectral sequences $\espse$.

\begin{prop}\label{prop:char-espse}
An object $(L, d_i, w_i,s_i)$ of $\lwb$ is in the essential image of $\Nn$ if and only if
$\Ker d_i = \Ker s_i$ for all $i\geq 0$.
\end{prop}

\begin{proof}
If $L$ is isomorphic to $\Nn X$, the description of the maps $w,s,d$ given in Proposition \ref{P:nerve_descr} proves that $\Ker d_i=\Ker s_i$ for all $i\geq 0$.

We have that $L$ is in the essential image of $\Nn$ if and only if $(\eta_L)_j$ is injective for all $j$, since it is surjective by Theorem~\ref{T:adj_QN}. \\
Let us assume that $\Ker d_i=\Ker s_i$ for all $i\geq 0$. We prove by induction on $r$ that $\Ker(\eta_L)_r=0$. Firstly note that  $\Ker(\eta_L)_0=0$ for all $L$ since $(\Nn(\Qq L))_0$ is isomorphic to $L_0$. Hence the result is valid for $r=0$. 
Let $r\geq 0$ and assume the statement holds for every $k\leq r$.
Let $a\in L_{r+1}$ be such that $(\eta_L)_{r+1}(a)=0$. Thus $(\eta_L)_r(w_{r+1}a)=0$ and
$(\eta_L)_r(w_{r+1}d_{r+1}a)=0$, and so $w_{r+1}a=0$ and $w_{r+1}d_{r+1}a=0$ since $(\eta_L)_r$ is injective. By assumption, the image of $a$ in $(\Qq L)_{r+1}$ is zero, thus there exist $\alpha,\beta$ such that $a=s_r\alpha+d_{r+1}s_r\beta$. Applying $w_{r+1}$ we get $0=w_{r+1}a=d_r\beta$.
So $\beta\in\Ker d_r\subseteq \Ker s_r$ and  hence $s_r\beta=0$. Applying $w_{r+1}d_{r+1}$ we get $0=w_{r+1}d_{r+1}a=d_r\alpha$, hence
 $\alpha\in\Ker d_r\subseteq \Ker s_r$  and so $s_r\alpha=0$. In conclusion $a=0$.
 \end{proof}

\begin{defi}
We denote by $\lwbe$ the full subcategory of $\lwb$ on objects satisfying $\Ker d_i = \Ker s_i$ for all $i\geq 0$.
\end{defi}

\begin{teo}\label{T:eqcat_espse}
The fixed point equivalence associated to the
adjunction $(\Qq, \Nn)$ gives an equivalence of categories $\lwbe\sim \espse$.
\end{teo}

\begin{proof}
It follows from Theorem~\ref{T:adj_QN}
and Proposition~\ref{prop:char-espse} that the adjoint functors $\Qq$ and $\Nn$ restrict to give an 
equivalence: 
\[\xymatrix{
 \lwbe\ar@<1ex>[rr]^-{\Qq} && \espse \ar@<1ex>[ll]^-{\Nn}_-{\sim} }.\qedhere \]
\end{proof}

We now note some properties of the category $\lwbe$.
In Proposition~\ref{P: charespse_rlp}, we explain how to characterize $\lwbe$ via a right lifting property. And in
Remark~\ref{local-lwbe}, we note that it can be viewed as a certain
subcategory of local objects.

\begin{prop}\label{P:Sacyclic}
\begin{enumerate}
\item For $L$ in $\lwbe$,
the $r$-bigraded complex $S_r(L)$ is the direct sum $s_{r-1}L_{r-1}\oplus d_r s_{r-1} L_{r-1}$ and $S_r(L)$ is acyclic for all $r$. 
\item For $L$ in $\lwbe$, the quotient map $L_r\to\Qq(L)_r$ is a quasi-isomorphism
for all $r$.
\item For $X$ an extended spectral sequence, the projection  $\rho_r\colon (\Nn X)_r\to X_r$ is a quasi-isomorphism for all $r$.
\item For $X$ an extended spectral sequence
satisfying $d_i=0$ for $0\leq i\leq r-1$ and $\varphi_i=1$ for $1\leq i\leq r$, the projection  $\rho_r\colon (\Nn X)_r\to X_r$ is an isomorphism.
\end{enumerate}
\end{prop}

\begin{proof}
Let $a\in L$ with $s_{r-1}a=d_rs_{r-1}b$. Then applying $w_r$ we get $d_{r-1}b=0$ hence $s_{r-1}b=0$ and so $s_{r-1}a=0$ too. So we see that
$S_r(L)$ is a direct sum. So it is the cone in $r$-bigraded complexes of the identity map of $s_{r-1}L_{r-1}$ and so acyclic.
Since the kernel of the quotient map $L_r\to \Qq(L)_r$ is $S_r(L)$, we deduce that the quotient map is a quasi-isomorphism for all $r$.
The projection $\Nn X\to X$ is the composition of the quotient map $\Nn X\to \Qq\Nn X$ followed by the counit $\epsilon_X$ of the adjunction $(\Qq,\Nn)$, which is an isomorphism. Since $\Nn X\in\lwbe$, the third statement follows.
The final statement can be seen directly from the description of $(\Nn X)_r$.
\end{proof}

\subsection{Limits and colimits in extended spectral sequences}

\begin{prop}\label{P:lwbe-lims-colims}
The category $\lwbe$ is a full reflective subcategory of $\lwb$. It has all small limits and colimits, with limits computed pagewise and colimits
via the reflector $\Nn\Qq$.
\end{prop}

\begin{proof}
By Theorem~\ref{T:adj_QN}, the counit of the adjunction $(\Qq, \Nn)$ is a natural isomorphism. 
It follows that the essential image of the
right adjoint functor $\Nn$ is a full reflective subcategory of $\lwb$, with reflector $\Nn\Qq$ left adjoint to the inclusion.
Proposition~\ref{prop:char-espse} identifies the essential image of $\Nn$ as $\lwbe$.

Note that $\lwb$ is bicomplete with small limits and colimits, calculated componentwise.
It follows that the inclusion creates limits and the full subcategory also has all colimits, obtained by applying the reflector to colimits in
$\lwb$~\cite[Proposition 4.5.15]{Riehl_context}.
\end{proof}

\begin{teo}\label{T:espse-lims-colims}
The category $\espse$ has all small limits and colimits, with
colimits computed pagewise. 
\end{teo}

\begin{proof}
The existence of small limits and colimits follows from the equivalence of categories of Theorem~\ref{T:eqcat_espse}.
In Proposition~\ref{P:espse_colimit}, we saw that colimits are pagewise. 
If $G\colon I\to\espse$  is a functor from a small category $I$, then $\lim_i G(i)=\Qq (\lim_i \Nn(G(i)))$.
\end{proof}

As we noted earlier, the category $\espse$ is not cotensored over $\modu$.
So, although it is complete, one should not expect $\Nn$ to have a right adjoint. Indeed, $\Nn$ does not preserve all colimits, but we do have the following result which will be useful for transfer of model structures later.

\begin{prop}\label{P:N_filteredcolim} The functor $\Nn\colon \espse\to\lwb$ preserves filtered colimits.
\end{prop}

\begin{proof} We recall that colimits are computed pagewise in $\espse$ and in $\lwb$.
Suppose that $i\colon \mathcal{C}\to \mathcal{C}'$ is the inclusion of a full subcategory, that $\colim_J (i\circ X)$
exists in $\mathcal{C}'$ and is (isomorphic to) an object
of $\mathcal{C}$. Then $\colim_J X$ exists and is given by the same object. Applying this to the inclusion of the essential image of $\Nn$ into $\lwb$, we
therefore need to show that for filtered colimits in $\lwbe$ the pagewise colimit lies in $\lwbe$.

Let $L\colon J\to\lwbe$ be a diagram in $\lwbe$ with $J$ a filtered category. The pagewise colimit
is
\[
	C:=\colim_J L(j)=(\oplus_{j\in J} L(j))/\sim
\]
where $a\sim a'$ for $a\in L(j), a'\in L(j')$ if there exist $f\colon j\to k, f'\colon j'\to k$ such that $L(f)(a)=L(f')(a')$.
Let $[a]=[(a_j)_j]\in C$, where finitely many $a_j$ are non-zero, say for $j\in J'$, a finite subset of $J$. 
Suppose $d_i[a]=0$ for some $i\geq 0$. Now since $J$ is filtered there exists $k\in J$ and $f_j\colon j\to k$ in $J$ for all $j\in J'$. 
Let $b= \Sigma_{j\in J'}L(f_j)(a_j)\in L(k)$. 
Then $a=\Sigma_{j\in J'} a_j\sim b$. Now $0=d_i[a]=d_i[b]=[d_ib]$. So $d_ib\sim 0$ and there exists $f\colon k\to k'\in J$ such that
$d_iL(f)(b)=L(f)(d_ib)=0$. Thus $L(f)(b)\in \Ker d_i=\Ker s_i$ in $L(k')$. And since $a\sim b\sim L(f)(b)$, we see that
$s_i[a]=0$. Thus $C$ satisfies $\Ker d_i=\Ker s_i$ for all $i\geq 0$ and so $C$ lies in $\lwbe$.
\end{proof}

\subsection{Characterizing spectral sequences as presheaves}\label{S:lwbs}

\begin{defi}
We denote by $\lwbs$ the subcategory of $\lwb$ given by the essential image under $\Nn$ of the category of spectral 
sequences $\spse$. 
\end{defi}

By restricting our adjunction we obtain an equivalence of categories  $\lwbs\sim \spse$. 
Equivalently,  $\lwbs$ is the full subcategory on those $L$ in $\lwbe$ such that $\Qq L$ is a spectral sequence.

The following diagram gives an overview of the relationships between these categories.
 \begin{center}  
    \begin{tikzcd}
\lwb \arrow[r, yshift=0.8ex, "\Qq"{name=Q}]&\espse \arrow[l, yshift=-0.6ex, "\Nn"{name=N}] \arrow[ld, yshift=-0.7ex, " "{name=D}]{}{}\\
\arrow[phantom, from=Q, to=N, , "\scriptscriptstyle\boldsymbol{\bot}"]
\lwbe\arrow[u, hook]\arrow[ru, yshift=0.8ex, " "{name=U}]{}{}&\\
\arrow[phantom, from=U, to=D, , "\scriptscriptstyle\boldsymbol{\sim}"]
\lwbs\arrow[u, hook]\arrow[r, yshift=0.7ex, " "{name=T}]{}{}&\spse\arrow[uu, hook] \arrow[l, yshift=-0.7ex, " "{name=B}]{}{}\\
       \arrow[phantom, from=T, to=B, , "\scriptscriptstyle\boldsymbol{\sim}"]                  
    \end{tikzcd}
      \end{center}
\vspace{-1cm}

\begin{prop}\label{P: charspseq}
Let $(L,d_i,w_i,s_i)$ be an object of $\lwbe$ and let $i\geq 0$.
\begin{enumerate}
\item The characteristic map $\varphi_{i+1}$ of $\Qq L$ is surjective if and only if
in $L$ we have $\Img w_{i+1}=\Ker d_i=\Ker s_i$.
\item The characteristic  map $\varphi_{i+1}$   of $\Qq L$  is injective if and only if
in $L$ we have $\Ker w_{i+1}=\Img s_i$.
\end{enumerate}
Thus an object $(L, d_i, w_i,s_i)$ of $\lwb$ lies in $\lwbs$ if and only if
$\Img w_{i+1}=\Ker d_i=\Ker s_i$ and $\Ker w_{i+1}=\Img s_i$
for all $i\geq 0$.
\end{prop}

\begin{proof} We may consider $L$ as $L=\Nn X$ with $X$ an extended spectral sequence. In particular $\Qq L$ is isomorphic to $X$ by Theorem~\ref{T:adj_QN}. We recall that $\Img w_{i+1}\subseteq \Ker s_i$ and $\Img s_i\subseteq \Ker w_{i+1}$.
Let $(x;y)$ in $(\Nn X)_i$. We have 
\begin{itemize}
\item $(x;y)\in \Img w_{i+1}$ if and only if $y=0$ (and thus $d_ix_i=0$) and there exists $x_{i+1}\in X_{i+1}$ such that $\varphi_{i+1}(x_{i+1})=[x_i]$.
\item $(x;y)\in \Ker s_i$ if and only if $y=0$ (and thus $d_ix_i=0$).
\end{itemize}
In consequence the characteristic map $\varphi_{i+1}$ of $\Qq L$ is surjective if and only if
in $L$ we have $\Img w_{i+1}=\Ker s_i=\Ker d_i$.

Let $(x;y)$ in $(\Nn X)_{i+1}$. We have 
\begin{itemize}
\item $(x;y)\in \Img s_{i}$ if and only if $x=0$ (and thus $x_{i+1}=0$). 
\item $(x;y)\in \Ker w_{i+1}$ if and only if $(x_0,\ldots,x_i)=(0,\ldots,0)$. 
\end{itemize}
In consequence $\Ker w_{i+1}=\Img s_i$ if and only if, for all $ x_{i+1}\in X_{i+1}$, $\varphi_{i+1}(x_{i+1})=0\Rightarrow x_{i+1}=0$.
So the characteristic map $\varphi_{i+1}$  of $\Qq L$  is injective if and only if
in $L$ we have $\Ker w_{i+1}=\Img s_i$.

The final statement follows since $L\in\lwbs$ if and only if
the characteristic map $\varphi_{i+1}$ of $\Qq L$ is an isomorphism for all $i\geq 0$.
\end{proof}

In Proposition~\ref{P:charspse_rlp}, we explain how to characterize $\lwbs$ via a right lifting property. 

\section{Truncation, shift and d\'ecalage}
\label{S:TSD}

In this section we will study shift and d\'ecalage functors.
Up to technicalities related to bigrading, the shift functor on
(extended) spectral sequences is easy to describe: it introduces a duplicate of the \(0\)-page at the start, but with trivial differential, and shifts all the other pages up one. It turns out that this functor has both a left and a right adjoint functor, giving two versions of d\'ecalage. The décalage functor \(\Dec\), right adjoint to the shift functor, is described similarly. For \(X\) an (extended) spectral sequence, up to bigrading, \(\Dec(X)_0\) is given by \(\Nn(X)_1\) and \(\Dec(X)_k\) by \(X_{k+1}\), for \(k\geq 1\).
We carry out a detailed study of how these functors arise from
structure on the underlying category of discs \(\D\).

To this end we introduce various truncated versions of the categories
$\D$, $\lwb$ and $\espse$ and set up functors and adjunctions relating these. 
Indeed all of these categories contain nested copies of themselves giving rise to
these adjunctions.

In the first part we consider truncated versions of $\D$, where we take
objects $(i,p,n)$ with $r\leq i\leq t$. We study
the embedding of this truncated version into $\D$, and the associated triple of adjoint functors at the level of linear witness books. 
The special case $r=t$ will be used for transfer of model category structures in Section~\ref{S:ModelsLWB}.

In the second part we
consider the case $t=\infty$. In this case, the embedding into $\D$ has both a left and a right adjoint. This triple of adjoint functors gives rise to a sequence of five adjoint functors on $\lwb$. Of these, we
 will see that to the leftmost three there is a corresponding triple of adjoints
  on $\espse$ and on $\spse$. 
  
  When $r=1$ we use the term d\'ecalage 
 for the functor on $\espse$ or $\spse$ corresponding to the embedding 
 and denote it by $\Dec$.  Its left adjoint is called the shift functor and the left adjoint of the shift is a variant of d\'ecalage, which we denote $\LDec$.
We explain our use of the term d\'ecalage by noting compatability with Deligne's notions for filtered complexes.
The presheaf approach gives important insight in determining the correct form of the d\'ecalage functor for
(extended) spectral sequences. 

This section is quite technical, involving many adjoint functors. These arise formally as left or right Kan extensions
as in Proposition~\ref{P:presheaf-yoga}. In several cases we give explicit descriptions of these functors on objects, omitting the corresponding descriptions for morphisms and omitting proofs which are routine. We defer an explicit description of the functor \(\LDec\) and its iterates to Appendix~\ref{App:adjoints}. This is the most delicate case
and we provide
 full details in Proposition~\ref{P:ldecr-shiftr-espse}.

\subsection{Truncation}
\label{SS:truncation}
We begin by considering truncated versions of the category $\D$ and the corresponding truncations of the category of linear witness books $\lwb$, where we have only pages $r$ to $t$, for some $0\leq r\leq t\leq \infty$. Truncation of $\D$ may be thought of as analogous to truncation of the simplex category, giving rise to (co)skeleta functors on simplicial sets.

\begin{notation}
    For $0\leq r\leq t\leq \infty$, we denote by $\tr_{[r,t]}\D$ 
the full subcategory of $\D$ generated by the objects $(i,p,n)$ for 
$r\leq i\leq t$. We write
$\tr_{\geq r}\D$ for $\tr_{[r,\infty]}\D$ and 
$\tr_{\leq r}\D$ for $\tr_{[0,r]}\D$ when convenient. Note that $\tr_{[r,r]}\D=\langle \delta_r\rangle$.

We use the corresponding notation  for $\lwb$. Indeed the category $\tr_{[r,t]}\lwb$ has for objects families of $i$-bigraded complexes $(L_i,d_i)$ with $r\leq i\leq t$ together with linear maps $w_{i+1}:L_{i+1}\to L_{i}$ of bidegree $(0,0)$ and $s_i:L_i\to L_{i+1}$ of bidegree $(1,1)$  for $r\leq i<t$ subject to the same relations as those for linear witness books. Note that we have $\tr_{[r,r]}\lwb=\widehat{\langle\delta_r\rangle}$.

We denote by $\iota_{[r,t]}\colon\tr_{[r,t]}\D\to \D$ the embedding.

\end{notation}

The functor $\iota_{[r,t]}$ induces adjunctions on the associated linear presheaf categories;
see Proposition~\ref{P:presheaf-yoga}.

\begin{defi}
We define 
    \[\Uu_{[r,t]}=\iota_{[r,t]}^*,\quad \Ff_{[r,t]}=(\iota_{[r,t]})_!\quad \text{ and }\quad \Rr_{[r,t]}=(\iota_{[r,t]})_*.\]
\end{defi}
\smallskip

So we have the adjunctions:
 \[
 \begin{tikzcd}
    \tr_{[r,t]}\lwb
\ar[rr, bend left, shift left=0.8, outer sep=1mm,  "\Ff_{[r,t]}",  "\perp"'] 
\ar[rr, xshift=-0.7ex, yshift=-0.5ex,
bend right, 
outer sep=1mm, pos=0.48, "\Rr_{[r,t]}"', "\perp"]
    & &\lwb 
\ar[ll, "\Uu_{[r,t]}" description] 
\end{tikzcd}
\]

We note that $\Uu_{[r,t]}:\lwb\to\tr_{[r,t]}\lwb$ is the forgetful functor.  
The following proposition describes the left adjoint functor $\Ff_{[r,t]}$ explicitly. For the interested reader, the description of the right adjoint
$\Rr_{[r,t]}$ is given in Proposition~\ref{P:RU-LWB}. The latter functor does not behave as nicely as the former functor with respect to $\lwbe$ 
and $\lwbs$. We omit proofs since the arguments are routine.

\begin{prop}
\label{P:FU-LWB}
Let $(L,d_i,w_i,s_i)$ be an object of $\tr_{[r,t]}\lwb$.
We have
 \[ \Ff_{[r,t]}(L)_i^{p,n}=\begin{cases}
 L_r^{p,n} & \text{ for } 0\leq i<r, \\
 L_i^{p,n}& \text{ for } r\leq i\leq t,\\
 \cone_i(L_t/B_t(L_t))^{p+t-i,n+t-i}& \text{ for } i>t, \end{cases}\]

with the maps $d_i^\Ff,w_i^\Ff,s_i^\Ff$ defined as
 \begin{align*}
     d_i^\Ff&=0& w_{i+1}^\Ff&=1 & s_i^\Ff&=0 & \text{ for } 0\leq i<r,\\
     d_i^\Ff&=d_i &  w_{i+1}^\Ff&=w_{i+1} & s_i^\Ff&=s_i & \text{ for } r\leq i< t, \\
     d_t^\Ff&=d_t &  w_{t+1}^\Ff&=\begin{pmatrix}0&d_t\end{pmatrix} & s_t^\Ff&=\begin{pmatrix}1\\0\end{pmatrix} & \text{ for } i=t, \\
     d_i^\Ff&=\begin{pmatrix}0&0\\1& 0\end{pmatrix} & w_{i+1}^\Ff&=\begin{pmatrix}0&0\\0& 1\end{pmatrix} & s_i^\Ff&=\begin{pmatrix}1&0\\0& 0\end{pmatrix} & \text{ for } i>t.
 \end{align*}

 The unit of the adjunction $(\Ff_{[r,t]},\Uu_{[r,t]})$ is the identity. The counit of the adjunction 
$\Ff_{[r,t]}\Uu_{[r,t]}(L)\rightarrow L$ takes the following form.
\begin{itemize}
\item For $0\leq i<r$ it corresponds to $w^{r-i}\colon L_r\rightarrow L_i$.
\item It is the identity at  page $i$ for $r\leq i\leq t$.
\item For $i>t$ it is the composite of the projection onto the first component $\cone_i(L_t/B_t(L_t))\to L_t/B_t(L_t)$  followed by $s^{i-t}\colon L_t/B_t(L_t)\to L_i$. \qed
\end{itemize}
\end{prop}
\smallskip

The translation functor of Definition~\ref{Def:shift_rbc} interacts
nicely with the truncations of the category $\lwb$.

\begin{prop}
    The translation functor applied to each page of a (truncated) linear witness book induces isomorphisms of categories
    \[\Tt\colon \tr_{[r,s]}\lwb\to\tr_{[r+1,s+1]}\lwb, \qquad \Tt\colon \lwb\to\tr_{\geq 1}\lwb\]
    and thus an isomorphism of categories
    \[\Tt^r\colon \lwb\to\tr_{\geq r}\lwb.\vspace{-0.75cm}\]
    \qed
\end{prop}
\medskip

\subsection{Special case of truncation}
Now we consider some extra properties that are satisfied by truncation $\iota_{[r,t]}\colon\tr_{[r,t]}\D\to \D$
in the case $t=\infty$. Recall that we write $\iota_{\geq r}$ for the embedding in this case.
This embedding functor
$\iota_{\geq r}\colon\tr_{\geq r}\D\to\D$ has both a left and a right adjoint, so that 
$\tr_{\geq r}\D$ is a reflective and coreflective subcategory of $\D$.
We therefore obtain via Proposition~\ref{P:presheaf-yoga} a chain of five adjoint functors at the level of the presheaf category $\lwb$ and its truncations. Most of our focus will be on the leftmost three of these five functors, since these are the ones that give rise to interesting functors on (extended) spectral sequences.

\begin{prop}\label{P:adjs_emb}
    The embedding $\iota_{\geq r}\colon\tr_{\geq r}\D\to\D$ admits a left and a right adjoint, denoted respectively 
    $\Ww_{\geq r}$ and 
    $\Ss_{\geq r}$ from $\D$ to $\tr_{\geq r}\D$. The left adjoint $\Ww_{\geq r}$
is given
on objects by
\begin{align*}
\Ww_{\geq r}(i,p,n)&=(r,p,n) \qquad\text{ if } i<r,\\
\Ww_{\geq r}(i,p,n)&=(i,p,n)\qquad\text{ if } i\geq r,
\end{align*}
and on generating morphisms by, for $i<r$,
\[
\Ww_{\geq r}(\omega_{i+1}^{p,n})= 1,
\quad
\Ww_{\geq r}(\delta_i^{p,n})= 0,
\quad
\Ww_{\geq r}(\sigma_i^{p,n})= 0,
\]
with $\Ww_{\geq r}$ acting as the identity on all other generating morphisms.

The right adjoint $\Ss_{\geq r}$
is given
on objects by
\begin{align*}
\Ss_{\geq r}(i,p,n)&=(r,p+r-i,n+r-i)&\text{ if } i<r,\\
\Ss_{\geq r}(i,p,n)&=(i,p,n)&\text{ if } i\geq r,
\end{align*}
and on generating morphisms by, for $i<r$,
\[
\Ss_{\geq r}(\omega_{i+1}^{p,n})= 0,
\quad
\Ss_{\geq r}(\delta_i^{p,n})= 0,
\quad
\Ss_{\geq r}(\sigma_i^{p,n})= 1,
\]
with $\Ss_{\geq r}$ acting as the identity on all other generating morphisms.

We have $\Ww_{\geq r}\iota_{\geq r}=\Ss_{\geq r}\iota_{\geq r}
=\Id_{\tr_{\geq r}\D}$.

The counit of the adjunction $(\Ww_{\geq r},\iota_{\geq r})$ is the identity. The unit of the adjunction
$\eta^{\Ww}_{(i,p,n)}: (i,p,n)\to \Ww_{\geq r}(i,p,n)$ has the following form.
\begin{align*}
\eta^{\Ww}_{(i,p,n)}&= \omega^{r-i}\colon (i,p,n)\to (r,p,n)& \text{ if } i<r,\\
    \eta^{\Ww}_{(i,p,n)}&=1 &\text{ if } i\geq r.
\end{align*}

The unit of the adjunction $(\iota_{\geq r},\Ss_{\geq r})$ is the identity. The counit of the adjunction
$\epsilon^{\Ss}_{(i,p,n)}\colon\Ss_{\geq r}(i,p,n)\to (i,p,n)$ has the following form.
\begin{align*}
\epsilon^{\Ss}_{(i,p,n)}&= \sigma^{r-i}\colon (r,p+r-i,n+r-i)\to (i,p,n)& \text{ if } i<r,\\
    \epsilon^{\Ss}_{(i,p,n)}&=1 &\text{ if } i\geq r.
\end{align*}
\end{prop}

\begin{proof}
    That the functors are well defined is immediate. To prove that $\Ww_{\geq r}$ is left adjoint to $\iota_{\geq r}$, we only need to show that
$\Ww_{\geq r}(\eta^\Ww_{(i,p,n)})$ and $\eta^\Ww_{\iota_{\geq r}(i,p,n)}$ are the identity maps, which is the case.
To prove that $\Ss_{\geq r}$ is right adjoint to $\iota_{\geq r}$, we only need to show that
$\Ss_{\geq r}(\epsilon^\Ss_{(i,p,n)})$ and $\epsilon^\Ss_{\iota_{\geq r}(i,p,n)}$ are the identity maps, which is the case.
\end{proof}

Thus we have a chain of five adjoint functors on (truncated) linear witness books
\[
\renewcommand*{\arraystretch}{1.3}
\begin{matrix}
(\Ww_{\geq r})_!&\dashv &\Ww_{\geq r}^*
&\dashv &\iota_{\geq r}^* &\dashv &\Ss_{\geq r}^*&\dashv &(\Ss_{\geq r})_*\\
&&\veq&&\veq&&\veq&&\\
&&\Ff_{\geq r}&&\Uu_{\geq r}&&\Rr_{\geq r}\\
\end{matrix}
\]

The explicit description of the functor $(\Ww_{\geq r})_!$ is given in Proposition~\ref{P:wexclamation}.

\subsection{Shift and d\'ecalage for (extended) spectral sequences}

Next we consider how these adjoint functors behave with respect to the subcategory $\lwbe$ and its truncations. We show that the leftmost three of the five adjoints restrict to this setting and thus give rise, via
the $(\Qq, \Nn)$ adjunction, to corresponding adjoint functors on $\espse$. We will show that they also restrict to $\spse$.

The translation functor $\Tt^r:\espse\to \tr_{\geq r}\espse$ is an isomorphism of categories. One can then define a pair of adjoint functors $(\Qq_{\geq r},\Nn_{\geq r})$ between the categories $\tr_{\geq r}\lwb$ and 
$\tr_{\geq r}\espse$, transferring the adjunction $(\Qq,\Nn)$ along the isomorphism $\Tt^r$.
In particular, the essential image of $\Nn_{\geq r}$ can be described as in Proposition ~\ref{prop:char-espse}; we denote this by $\tr_{\geq r}\lwbe$. This category has for objects the truncated linear witness books satisfying $\Ker s_i=\Ker d_i$, for all $i\geq r$.

\begin{prop}
    The functors $\Uu_{\geq r}$ and $\Ff_{\geq r}$ restrict to the categories $\lwbe$ and $\tr_{\geq r}\lwbe$ and to the categories $\lwbs$ and $\tr_{\geq r}\lwbs$.
\end{prop}

\begin{proof}
   The functor $\Uu_{\geq r}:\lwb\to\tr_{\geq r}\lwb$ is the forgetful functor, hence it sends $\lwbe$ to $\tr_{\geq r}\lwbe$ and $\lwbs$ to $\tr_{\geq r}\lwbs$. Let $L$ be a truncated linear witness book. For $0\leq i<r$ we have, in $\Ff_{\geq r} L$, $\Ker d_i=\Ker s_i=\Img w_{i+1}=L_r$
   and $\Img s_i=\Ker w_{i+1}=0$. For $i\geq r$ we have, in $\Ff_{\geq r} L$,  all the structural morphisms 
   coincide with those in $L$. In conclusion, if $L\in\tr_{\geq r}\lwbe$ or in $\tr_{\geq r}\lwbs$ then $\Ff_{\geq r} L$ is in $\lwbe$ or in $\lwbs$, by Propositions \ref{prop:char-espse} and \ref{P: charspseq}.
\end{proof}

The corresponding statement for the functor $(\Ww_{\geq r})_!$ can be found in Proposition~\ref{P:wwrestrict}.

Next we introduce notation for the triple of adjoint functors we are interested in on
$\espse$. These correspond to the functors
$(\Ww_{\geq r})_!, \Ff_{\geq r}$ and $ \Uu_{\geq r}$, transferred via the $(\Qq, \Nn)$ adjunction and composed with suitable translations in order to give functors $\espse\to\espse$.

\begin{notation}
    We introduce the following functors
    $\espse\to\espse$:
    \begin{align*}
        \LDec&=\Tt^{-1}\circ \Qq_{\geq 1}\circ (\Ww_{\geq 1})_!\circ \Nn,\\
        \Shift&=Q\circ \Ff_{\geq 1}\circ \Nn_{\geq 1}\circ \Tt,\\
        \Dec&=\Tt^{-1}\circ \Qq_{\geq 1}\circ \Uu_{\geq 1}\circ \Nn.
    \end{align*}
\end{notation}

We denote by $\tr_{\geq r}\espse$ the category whose objects are collections of $i$-bigraded complexes $(X_i,d_i)$ for $i\geq r$ endowed with characteristic maps $\varphi_{i+1}\colon X_{i+1}\to H(X_i)$ for $i\geq r$. And $\tr_{\geq r}\spse$ denotes the corresponding truncation of $\spse$.

\begin{teo}\label{T:ldec-shift-dec}
    On the category $\espse$ we have adjunctions
    \[\LDec\dashv \Shift\dashv \Dec.\]
\end{teo}

\begin{proof}
    We give the details for the $\Shift\dashv \Dec$ adjunction.
   Using that the unit of the adjunction $L\to \Nn_{\geq 1}\Qq_{\geq 1}(L)$ is an 
   isomorphism for $L\in \tr_{\geq 1}\lwbe$, that $\Uu_{\geq 1}$ sends $\lwbe$ to $\tr_{\geq 1}\lwbe$ and that $\Nn_{\geq 1}$ is fully faithful, we obtain the following sequence of bijections:
      \begin{align*}
       \Hom_{\espse}(\Shift(X),Y)&\cong
       \Hom_{\lwb}(\Ff_{\geq 1}\Nn_{\geq 1} \Tt(X), \Nn(Y))\\
       &\cong\Hom_{\tr_{\geq 1}\lwb}(\Nn_{\geq 1}\Tt (X),\Uu_{\geq 1}\Nn (Y))\\
       &\cong \Hom_{\tr_{\geq 1}\lwb}(\Nn_{\geq 1}\Tt(X),\Nn_{\geq 1}\Qq_{\geq 1}\Uu_{\geq 1}\Nn (Y))\\
       &\cong \Hom_{\tr_{\geq 1}\espse}(\Tt (X),\Qq_{\geq 1}\Uu_{\geq 1}\Nn (Y))\\
       &\cong \Hom_{\tr_{\geq 1}\espse}(\Tt(X),\Tt\Dec(Y))\\
       &\cong \Hom_{\espse}(X, \Dec(Y))
   \end{align*}
   The same arguments prove the other adjunction, using that $\Ff_{\geq 1}$ sends $\tr_{\geq 1}\lwbe$ to $\lwbe$.
\end{proof}

We write $F^r$ for the $r$-fold iterated composite $F\circ F \circ F \dots \circ F$ of a functor $F$.
The following proposition gives explicit descriptions of the iterated composites $\Shift^r$
and $\Dec^r$.

\begin{prop}\label{P:shiftr-decr-espse}
For $X$ a truncated extended spectral sequence in $\tr_{\geq r}\espse$, we have
    \[(\Shift^r \Tt^{-r} X)_i=\begin{cases}
(X_r,0) & \text{ for }  0\leq i< r,\\
(X_i,d_i) & \text{ for } r\leq i.\\
    \end{cases}\]
    The characteristic  maps $\varphi_i$ are  the identity morphism for $1\leq i\leq r$ and coincide with $\varphi_i^X$ for $r+1\leq i$. If $X$ is in $\spse$ so is $\Shift^r( X)$.
    
    \end{prop}

    \begin{proof}
         Obtaining the descriptions of the functor is straightforward. The statement about restriction to
         $\spse$ follow from the descriptions of the characteristic maps.
         \end{proof}

\begin{teo}\label{T:shiftdectheorem}
 On the category $\espse$, and on the category $\spse$, we have adjunctions
    \[\LDec^r\dashv \Shift^r\dashv \Dec^r,\]
with isomorphisms:
   \begin{align*}
        \LDec^r&\cong\Tt^{-r}\circ \Qq_{\geq r}\circ (\Ww_{\geq r})_!\circ \Nn,\\
        \Shift^r&\cong Q\circ \Ff_{\geq r}\circ \Nn_{\geq r}\circ \Tt^r,\\
        \Dec^r&\cong \Tt^{-r}\circ \Qq_{\geq r}\circ \Uu_{\geq r}\circ \Nn.
    \end{align*}
\end{teo}

\begin{proof}
    The first part is immediate from Theorem~\ref{T:ldec-shift-dec}, together with the fact that the functors restrict to
    $\spse$. The expression for $\Shift^r$ as a composite is easy to check directly with the explicit description of Proposition~\ref{P:shiftr-decr-espse}. Then the composite expressions for $\LDec^r$ and
    $\Dec^r$ give respectively left and right adjoints to $\Shift^r$ by the same argument as in the proof of
    Theorem~\ref{T:ldec-shift-dec}, thus establishing the other isomorphisms. 
\end{proof}

The following gives an explicit description of the functor 
\(\Dec^r\).

\begin{cor}\label{C:decr-espse}
    For $X$ an extended spectral sequence in $\espse$, we have in $\tr_{\geq r}\espse$
     \[(\Tt^r\Dec^rX)_i=\begin{cases}
(\Nn X)_r & \text{ for }  i=r,\\
(X_i,d_i) & \text{ for } r< i.\\
    \end{cases}\]
   The characteristic maps $\varphi_i$ coincide with $\varphi_i^X$ for $r+1<i$ and $\varphi_{r+1}$ is the following composite
 \[\begin{tikzcd}
X_{r+1}\arrow[r,"\varphi^X_{r+1}"]& H(X_r)\arrow[r,"H(\rho_r)^{-1}"] & H((\Nn X)_r).
   \end{tikzcd}\] 
   If $X$ is in $\spse$ then so is $\Dec^r(X)$.
   
The unit of the adjunction
$(\Shift^r, \Dec^r)$ is an isomorphism.
    \end{cor}

    \begin{proof}
         By Theorem \ref{T:shiftdectheorem}, we have $\Tt^r\Dec^rX=\Qq_{\geq r}\Uu_{\geq r}\Nn X$ and the description follows.
         The statement about restriction to
         $\spse$ follows from the description of the characteristic maps.
         That the unit $\eta$ is an isomorphism follows from
         part (4) of Proposition~\ref{P:Sacyclic}, with $(\eta_X)_0=(\rho_r)^{-1}$ and $(\eta_X)_i=1$ for $i>0$.
    \end{proof}

The corresponding description for the functor $\LDec^r$ can be found in Proposition~\ref{P:ldecr-shiftr-espse}.

\subsection{Combatibility with shift and d\'ecalage for filtered complexes}\label{S:fcx-dec}

In this section we look at compatibility of the functors 
 $\Shift$ and $\Dec$ with the functors 
 $S$ and $\Dec$ on the category of filtered complexes $\fcpx$. These were introduced by Deligne 
in~\cite[Section 1.3]{DeHII}; a useful exposition of their properties is given in
~\cite[Section 2.3]{CG2}.  We write $E$ for the functor associating a spectral sequence to a filtered complex
and we will use the notation of
~\cite[Section 2]{CELW20} for constructions related to the spectral sequence of a filtered complex.

Deligne also studied the \emph{left} adjoint  $\Dec^*$ of
$S$. Compatibility of 
$\LDec$ and $\Dec^*$ is established in Proposition~\ref{prop: ldec-fcx}.

\begin{prop}
We have a commutative diagram
\[\xymatrix{
\fcpx \ar[r]^-{E}\ar[d]_{S}&\spse \ar[d]^{\Shift}\\
\fcpx \ar[r]_-{E}& \spse  }\]
\end{prop}

\begin{proof}
This is a matter of direct calculation with the definitions. Let $C$ be a filtered complex. For $i\geq 1$, we have $E(SC)_i=\Tt(E(C)_{i-1})$ and for
$i=0$ we obtain $E(SC)_0=\Tt(E(C)_{0})$ with zero differential.
\end{proof}

\begin{prop}
Let $C$ be a filtered complex.
There is a natural isomorphism of spectral sequences $E(\Dec(C))\to \Dec(E(C))$.
\end{prop}

\begin{proof}
Firstly consider the functor $\Tt^{-1}\circ \tr_{\geq 1}\colon \spse\to\spse$ 
which just shifts down the pages, differing from $\Dec$ only by information related to the $0$-page.
From~\cite[Proposition 1.3.4]{DeHII},
we have a canonical morphism of spectral sequences 
\[u\colon E(\Dec(C))\to \Tt^{-1}\circ \tr_{\geq 1}(E(C))\] which is a surjection on the $0$-page and an isomorphism from the $1$-page onwards.
In fact $u_0\colon  E(\Dec(C))_0\to \Tt^{-1}\circ \tr_{\geq 1}(E(C))_0$ factors via  $\Dec(E(C))_0$
and we use this to produce the required isomorphism. 
Recall that 
\[
E(\Dec(C))_0^{p,p+n}=\Tt^{-1}(Z_1^{p,p+n}(C)/Z_1^{p-1,p+n-1}(C))
\]
with notation as in~\cite[Section 2]{CELW20}. 
Define 
\[\tilde{u}_0\colon  E(\Dec(C))_0\to
\Dec(E(C))_0=\Tt^{-1}(\Nn E(C))_1
\]
 by
$[x]\mapsto ([x]_0, [x]_1; [dx]_0, [dx]_1)$.
Direct checks show that $\tilde{u}_0$ is a well-defined isomorphism of $0$-bigraded complexes. Furthermore 
$u_0=\rho_1\tilde{u}_0$ and
we have the following commutative diagram.
 \[
 \begin{tikzcd}
     \Tt^{-1}H(E(C)_1)&&&&&\Tt^{-1}E(C)_{2} \arrow[lllll, "\Tt^{-1}\varphi_{2}^{E(C)}"']
     \arrow[dlllll, "\hspace{-1cm}\Tt^{-1}(H(\rho_1)^{-1}\varphi_2^{E(C)})"]\\
     \Tt^{-1}H(\Nn(E( C))_1)\arrow[u, "\Tt^{-1}H(\rho_1)"]\\
     H( E(\Dec(C))_0)\arrow[u, "H(\tilde{u}_0)"]&&&&&E(\Dec(C))_1
     \arrow[uu, "{u_1}"']
     \arrow[lllll, "\varphi_1^{\Dec(C)}"']
 \end{tikzcd}
 \]

Thus $\tilde{u}=(\tilde{u}_0, u_1, u_2, \dots): E(\Dec(C))\to \Dec(E(C))$ is an isomorphism of spectral sequences.
\end{proof}

This compatibility with d\'ecalage for filtered complexes clarifies
the relation of our construction to other uses of the term d\'ecalage in connection with spectral sequences, such as in~\cite{BHS, Lev15}.
Sometimes the term d\'ecalage of a spectral sequence is used to refer to 
d\'ecalage of some underlying filtered object. Other times it refers to the
corresponding process of \emph{turning the page} of the spectral sequence, up to reindexing, generally only specifying the resulting spectral sequence from the $1$-page onwards. Here we have provided the two versions of d\'ecalage as  functors on the category
of spectral sequences, including the (not completely obvious) descriptions of the $0$-page of each.

\section{Model category structures on extended spectral sequences}\label{S:ModelsESpSe}

In this section we establish two different families of model category structures on $\espse$. 
In the first part, we work with model category structures on $\lwb$. The main result here is Theorem~\ref{T:trans_model_lwb},
establishing such a structure $\lwb_r$ for each $r\geq 0$. 
Then we move to extended spectral sequences and we first use right transfer from  the model category structure $\lwb_r$, resulting in the model structure
$\espse_r$ of Theorem~\ref{T:modelespse2}, which we show is Quillen equivalent to the projective model structure on $0$-bigraded complexes. Then in Theorem~\ref{T:modelespse3} we give a different model category structure on $\espse$, denoted $\espse_r'$, where weak equivalences $\Ee'_r\subset \Ee_r$ are those maps that are quasi-isomorphisms at page $r$ and isomorphisms at pages $i>r$. The advantage of this class of morphisms is that it  more accurately reflects the weak equivalences in spectral sequences.
Also we have a hierarchy of weak equivalences: $\Ee_r'\subset \Ee_{r+1}'$ for all $r$.

\subsection{Model category structures from transfer on linear witness books}
\label{S:ModelsLWB}
For each $r\geq 0$, we present  a model category structure on the category $\lwb$ of linear witness books, closely related to the projective model structure
on $r$-bigraded complexes, utilising a criterion for right transfer due to Drummond-Cole and Hackney ~\cite{DCH19}. Such a model structure can be directly obtained by
a right transfer and then we explain how to modify the 
transferred structure to give the model structure of Theorem~\ref{T:trans_model_lwb}.

We consider the special case $r=t$ of the adjunctions of Section~\ref{SS:truncation} for linear witness books
and abbreviate to 
 $\Ff_{[r]}$, $\Uu_{[r]}$ the functors $\Ff_{[r,r]}$, $\Uu_{[r,r]}$. 
Recalling that  $\tr_{[r,r]}\lwb=\widehat{\langle\delta_r\rangle}$,
we focus on the adjunction:
 \[
 \begin{tikzcd}
 \widehat{\delrcat}\ar[rr, shift left=0.8, yshift=0.8ex, "\Ff_{[r]}", "\perp"' ] 
 & & \lwb \ar[ll, yshift=-0.8ex, "\Uu_{[r]}"] 
\end{tikzcd}
\]

The following is immediate from the description of $\Ff_{[r]}$, see Proposition~\ref{P:FU-LWB}. 

\begin{prop}\label{P:homologyF} 
Let $C$ be an $r$-bigraded complex.
The homology of the linear witness book $\Ff_{[r]}(C)$ is given by
\[
H((\Ff_{[r]}(C))_m)=\begin{cases}
C&\text{ if $0\leq m\leq r-1$},\\
H(C) &\text{ if $m=r$},\\
0 &\text{ if $m> r$}.
\end{cases}\vspace{-0.4cm}
\]
\qed
\end{prop}

\begin{prop}\label{P:Finspse} 
 Let $C$ be an $r$-bigraded complex. The linear witness book $\Ff_{[r]}(C)$ is an object of $\lwbe$ if and only if $H(C)=0$ if and only if it is an object of $\lwbs$.
\end{prop}

\begin{proof} Note that for every $i\not=r$ we have $\Ker s_i=\Ker d_i$. For $i=r$, we have $\Ker s_r=B_rC$ and $\Ker d_r=Z_rC$. Hence the first equivalence follows from 
Proposition \ref{prop:char-espse}. In addition, we have $\Img w_{i+1}=\Ker s_i$ for all $i$ and $\Img s_i=\Ker w_{i+1}$ for all $i\not=r$. We have 
$(s_r(C))^{pq}=(C/B_rC)^{p-1,q-1}$ and 
$(\Ker w_{r+1})^{pq}=(C/B_rC)^{p-1,q-1}\oplus (Z_rC/B_rC)^{p+r,q+r-1}$. Thus the second equivalence follows from Proposition \ref{P: charspseq}.
\end{proof}

We write $\Zz(r,p,n)\in\lwb$ for the object
representing $\Ker d_r$; see~\ref{N:rep_ker} 
for more details. Since
$\delta_r^2=0$ we have a natural map $\delta_r\colon \Zz(r,p,n)\rightarrow \Yy(r,p+r,n+1-r)$.
Note that because $\Ff_{[r]}$ preserves colimits, we have $\Zz(r,p,n)=\Ff_{[r]}(\kk^{p,n})$. In particular Proposition~\ref{P:Finspse} shows that
$\Zz(r,p,n)$ does not lie 
in $\lwbe$; see also Lemma~\ref{L:Z_notin_lwbe}.

\begin{notation}\label{N:IrJr}
We introduce notation for various sets of morphisms of $\lwb$.
\begin{align*}
I_r&\colon =\left\{\delta_r\colon \Zz(r,p,n)\rightarrow \Yy(r,p+r,n+1-r)\right\}_{p,n\in\ZZ},\\
J_r&\colon =\left\{0\lra \Yy(r,p,n)\right\}_{p,n\in\ZZ},\\
I_{\leq r}&\colon =I_r\cup \cup_{i=0}^{r-1} J_i,\\
J_{\leq r}&\colon =\cup_{i=0}^r J_i.
\end{align*}
\end{notation}

Recall the model structure on $r$-bigraded complexes of Proposition~\ref{P:modcatrbc}, with generating sets $I, J$ and
note that we have $\Ff_{[r]}(I)=I_r$ and $\Ff_{[r]}(J)=J_r$.

\begin{defi} Let $K,L \in\lwb$. A morphism $f\colon K\rightarrow L$ is an {\sl $r$-fibration} if $f_i$ is bidegreewise surjective
for $0\leq i\leq r$. It is an {\sl $r$-weak equivalence} if $H(f_r)$ is an isomorphism. We denote by $\Fib_r$ the class of $r$-fibrations and by $\Ee_r$ that of $r$-weak equivalences.
\end{defi}

\begin{teo}\label{T:trans_model_lwb}
There is a cofibrantly generated model structure with generating sets $I_{\leq r}, J_{\leq r}$ on $\lwb$, denoted $\lwb_{r}$. This has fibrations 
$\Fib_r$ and weak equivalences $\Ee_r$.
\end{teo}

\begin{proof} 
We first show that there is a model category structure cofibrantly generated by the sets $I_r$ and $J_r$.  For this, note that
the functor $\Uu_{[r]}$ admits a left adjoint  $\Ff_{[r]}$ and a right adjoint $\Rr_{[r]}$ (described in Proposition~\ref{P:RU-LWB}). Moreover  $\Uu_{[r]}\Ff_{[r]}$ and $\Uu_{[r]}\Rr_{[r]}$ are both the identity functor and hence form a Quillen adjunction. By ~\cite[Theorem 2.3]{DCH19}, there is a cofibrantly generated model structure, where the generating classes of cofibrations and acyclic cofibrations are $I_r=\Ff_{[r]}(I)$ and $J_r=\Ff_{[r]}(J)$ respectively. 
In this model structure, the class of weak equivalences is 
$\Ee_r$ and the class of fibrations consists of those maps $f\colon K\to L$ such that $f_r$ is surjective.

Now we modify this model structure. 
We claim that $I_{\leq r}$ and $J_{\leq r}$ are the generating cofibrations and generating acyclic cofibrations
 respectively of a model category structure as in the statement. We follow the treatment in~\cite{Hovey}. 
 A morphism has the right lifting property with respect to $J_{\leq r}$ if and only if
$f_i$ is surjective for $0\leq i\leq r$.  It is clear that $\Ee_r$ satisfies the two-out-of-three property
 and is closed under retracts. The smallness conditions are satisfied. 
 We have $\cl{I_{\leq r}}{inj}=\Ee_r\cap \cl{J_{\leq r}}{inj}$ just as in the proof of~\cite[Theorem 3.16]{CELW20}.
It remains to prove that the class of maps having the left lifting property with respect to $r$-fibrations is in $\Ee_r$. 
The technique we use here is that of~\cite{Fausk}. Let $f\colon L\to M$ be such a map. Consider the diagram
\[\xymatrix{ L\ar[r]^-{\vmat{1}{0}}\ar[d]_f & L\oplus M'\ar[d]^{(f,p)}\\ M\ar[r]^{=} & M}\]
where $M'=\oplus_{i=0}^r \Ff_{[i]}\cone_i\Uu_{[i]}M$ (see Definition~\ref{D:rcone1} for the cone)  and $p=p_0+\ldots+p_r$ with 
\[
	p_i\colon\Ff_{[i]}\cone_i\Uu_{[i]}M \rightarrow \Ff_{[i]}\Uu_{[i]}M \rightarrow M
\]
 induced by the projection $\pi_i$ in 
$\ibc$ followed by the counit of the adjunction $( \Ff_{[i]}, \Uu_{[i]})$.
From Proposition~\ref{P:homologyF} and  the fact that in $r$-bigraded complexes the $r$-cone is acyclic we have
\[
H((\Ff_{[i]}\cone_i\Uu_{[i]}M)_m)=\begin{cases}
\cone_i(M_i) & \text{if $0\leq m\leq i-1$,}\\
0 & \text{otherwise}.
\end{cases}
\]
In particular the $r$-bigraded complex $M'_r$ is acyclic.

From the descriptions of the functor $\Ff_{[i]}$ and the counit of the adjunction
in Proposition~\ref{P:FU-LWB} we have that 
$(\Ff_{[i]}\Uu_{[i]}M)_i\rightarrow M_i$ is the identity. Hence the $i$-page of $p_i$ is the projection from the $i$-cone of $\Uu_{[i]}M$ 
to $\Uu_{[i]}M$.
It is bidegreewise surjective. It follows that  
$p_0+\dots + p_r $ is an $r$-fibration.
Thus $(f,p)$ is an $r$-fibration and a lift exists in the diagram, 
say $(g,h)\colon M\rightarrow L\oplus M'$.
Considering the homology on the $r$-page in the diagram we see
that $H(f_r)$ is an isomorphism with inverse $H(g_r)$.
\end{proof}

\begin{teo}\label{T:eq_lwbr_rbc} The following adjunction is a Quillen equivalence.
\[\xymatrix{
 \widehat{\delrcat}\ar@<1ex>[rr]^-{\Ff_{[r]}} && \lwb_r \ar@<1ex>[ll]^-{\Uu_{[r]}}_-{\perp} }
 \]
In addition, the adjunction $\Ff_{\geq 1}\circ \Tt\dashv\Tt^{-1}\circ  \Uu_{\geq 1}$  gives a Quillen equivalence between 
$\lwb_{r}$ and  $\lwb_{r+1}$ for all
$r\geq 0$.
\end{teo}

\begin{proof} 
It is clear that $\Uu_{[r]}$ preserves fibrations and acyclic fibrations, so it is a Quillen adjunction and we have seen in Proposition~\ref{P:FU-LWB} that the unit of the adjunction is the identity. 
Furthermore the right adjoint $\Uu_{[r]}$ creates weak equivalences, so this is a Quillen equivalence.

For $\Ff_{\geq 1}\circ \Tt\dashv\Tt^{-1}\circ  \Uu_{\geq 1}$, we have that $\Tt^{-1}\circ  \Uu_{\geq 1}: \lwb_{r+1}\to\lwb_r$ 
preserves fibrations and weak equivalences. 
The unit is a natural isomorphism and 
 the counit of the adjunction
$\epsilon\colon \Ff_{\geq 1}\circ  \Uu_{\geq 1}\to \Id_{\lwb}$, given
by $\epsilon_L=(w^L_1, 1, 1, \dots)$,
 is in $\Ee_i$ for all $i\geq 1$.
\end{proof}

So all the model structures on $\lwb$ in this section are Quillen equivalent to the projective model structure on $0$-bigraded complexes.

\subsection{Model category structures from transfer on extended spectral sequences}\label{S:model-espse}
We will use right transfer to obtain
a model category structure on $\espse$, denoted $\espse_r$,
from the model category $\lwb_r$ described in Theorem~\ref{T:trans_model_lwb}.

\begin{teo}\label{T:modelespse2}  There is a cofibrantly generated model category structure on $\espse$ denoted $\espse_{r}$ with fibrations
$\Fib_r$, the class of maps $f:X\to Y$ such that $f_i$ is bidegreewise surjective for $0\leq i\leq r$, and weak equivalences $\Ee_r$, the class of maps $f:X\to Y$ such that $f_r$ is a quasi-isomorphism.
This model category structure is Quillen equivalent to $\lwb_r$.
\end{teo}

\begin{proof}
We use right transfer of the cofibrantly generated model category structure $\lwb_r$ by the adjunction
\[\xymatrix{
 \lwb\ar@<1ex>[rr]^-{\Qq} && \espse \ar@<1ex>[ll]^-{\Nn}_-{\perp} }\]
The transferred weak equivalences are maps $f\colon X\to Y$ such that $N(f)\in\Ee_r$, which is equivalent to $f\in\Ee_r$
by Proposition~\ref{P:Sacyclic}. The
fibrations are those maps $f\colon X\to Y\in \espse$ such that $\Nn(f)_i$ is surjective for $0\leq i\leq r$, or equivalently $f_i$ is surjective for $0\leq i\leq r$ by Proposition~\ref{P:Nsurj}. 
By Proposition~\ref{P:N_filteredcolim}, $\Nn$ preserves filtered colimits.
It then remains to show that any map $f:X\to Y$ having the left lifting property with respect to maps in $\Fib_r$ is in $\Ee_r$.
Consider the following diagram in $\espse$
\[\xymatrix{ X\ar[r]^-{\vmat{1}{0}}\ar[d]_f & X\oplus Y'\ar[d]^{(f,p)}\\ Y\ar[r]^{=} & Y}\]
where $Y'$ is the following extended spectral sequence:
\[Y'_m=\begin{cases} Y_m&\text{ if } m<r,\\
\cone_r(Y_r)& \text{ if } m=r,\\
0& \text{ if } m>r,\end{cases}\]
with $\varphi_{r}^{Y'}=\varphi_r^Y\circ\pi_r$ where $\pi_r\colon \cone_r(Y_r)\to Y_r$ is the projection of Definition \ref{D:rcone1}, and where $p$ is the identity on pages $i<r$, it is $\pi_r$ on page $r$ and then it is zero. So $(f,p)\in\Fib_r$ and a lift exists in the diagram, proving that $H(f_r)$ is an isomorphism.
Note that, because $Y'_m=0$ for $m>r$, the existence of the lift also implies that $f_m$ is an isomorphism for $m>r$.
In conclusion, we have a sequence of Quillen adjunctions
\[\xymatrix{
 \widehat{\delrcat}\ar@<1ex>[rr]^-{\Ff_{[r]}} && \lwb_r \ar@<1ex>[ll]^-{\Uu_{[r]}}_-{\perp} 
 \ar@<1ex>[rr]^-{\Qq} && \espse_r  \ar@<1ex>[ll]^-{\Nn}_-{\perp}
 }
 \]
 where the left adjunction is a Quillen equivalence. By the two-out-of-three property, proving that the right adjunction is a Quillen equivalence is equivalent to proving that the adjunction $(\Kk_r,\Nn_r)$ with $\Kk_r=\Qq\circ\Ff_{[r]}$ and $\Nn_r=\Uu_{[r]}\circ\Nn$
\[\xymatrix{
 \widehat{\delrcat}\ar@<1ex>[rr]^-{\Kk_r} && \espse_r  \ar@<1ex>[ll]^-{\Nn_r}_-{\perp}
 }
 \]
 is a Quillen equivalence.
 Explicitly, for every $r$-bigraded complex $C$, we have
 $\Kk_r(C)_i=C$ if $0\leq i\leq r$ and $0$ if $i>r$. The differential on each page is $0$ except at page $r$ where it coincides with that of $C$. The characteristic maps are either $1_C$ or $0$.
 By the definition of weak equivalences in $\espse_r$ we have that the composite $\Nn_r$ creates weak equivalences. Hence it is enough to prove that the unit of the adjunction $C\to \Nn_r(K_r C)$ is a quasi-isomorphism. Indeed it is an isomorphism.
\end{proof}

\begin{rmk}\label{R:newequiv}
    The proof of the above theorem shows that any map $f$ having the left lifting property with respect to $r$-fibrations, in addition to satisfying that $f_r$ is a quasi-isomorphism, also satisfies that $f_i$ is an isomorphism for $i>r$.
\end{rmk}

\begin{cor}
    If $L$  is a cofibrant linear witness book in $\lwb_r$, then the projection $L_r\to \Qq(L)_r=L_r/S_r(L)$ is a quasi-isomorphism.
\end{cor}

\begin{proof}
    Since the adjunction $(\Qq,\Nn)$ is a Quillen equivalence, for every cofibrant linear witness book $L$ in $\lwb_r$ and for any (fibrant) extended spectral sequence $X$ if $\alpha: \Qq(L)\to X$ is a weak equivalence in $\espse_r$, then 
    $\Nn(\alpha)\circ \eta_L$ is a weak equivalence in $\lwb_r$. Since every object $X$ is fibrant in $\espse_r$, one can pick $\alpha$ to be the identity of $\Qq(L)$. Thus $\eta_L$ is a weak equivalence. The projection $L_r\to \Qq(L)_r=L_r/S_r(L)$ is the $r$-page of the composite of $\eta_L$ with the projection $\Nn\Qq L\to\Qq(L)$ which is a quasi-isomorphism by Proposition~\ref{P:Sacyclic}, thus a quasi-isomorphism. 
\end{proof}

The next corollary is a r\'esum\'e of the model category structures we
have encountered so far and the Quillen equivalences between them.

\begin{cor}\label{C:w-dec-espse} 
  In the diagram of model categories below, every adjunction is a Quillen equivalence.
 
 \begin{center}  
    \begin{tikzcd}[row sep=10ex]
\espse_0 \arrow[rr, "\Shift^r" description]  
\arrow[d, xshift=0.9ex, "\Nn"{name=N}]
&&\espse_r \arrow[ll, yshift=2.5ex, "\LDec^r"'{name=LDec}, "\scriptscriptstyle\boldsymbol{\bot}"]
\arrow[ll, yshift=-2.5ex, "\Dec^r"{name=Dec}, "\scriptscriptstyle\boldsymbol{\bot}"']  
\arrow[rr, "\Shift" description] 
\arrow[d, xshift=0.9ex, "\Nn"{name=N2}]
&&\espse_{r+1} 
\arrow[ll, yshift=2.5ex, "\LDec"'{name=LDec}, "\scriptscriptstyle\boldsymbol{\bot}"]
\arrow[ll, yshift=-2.5ex, "\Dec"{name=Decd}, "\scriptscriptstyle\boldsymbol{\bot}"'] 
\arrow[d, xshift=0.9ex, "\Nn"{name=N3}] \\
\lwb_{0} \arrow[u, xshift=-0.9ex, "\Qq"{name=Q}] 
\arrow[rr, yshift=0.9ex, "(\Ff_{\geq 1}\circ \Tt)^r"] 
\arrow[d, xshift=0.9ex, "\Uu_{[0]}"{name=P}]
&&\lwb_r\arrow[u, xshift=-0.9ex, "\Qq"{name=Q2}]  
\arrow[ll, yshift=-1.25ex, "(\Tt^{-1}\circ  \Uu_{\geq 1})^r"{name=Dec2}, "\scriptscriptstyle\boldsymbol{\bot}"'] 
\arrow[rr, yshift=0.9ex, "\Ff_{\geq 1}\circ \Tt"]  
\arrow[d, xshift=0.9ex, "\Uu_{[r]}"{name=P2}]&&\lwb_{r+1}
\arrow[u, xshift=-0.9ex, "\Qq"{name=Q3}]  
\arrow[ll, yshift=-1.25ex, "\Tt^{-1}\circ  \Uu_{\geq 1}"{name=Dec2d}, "\scriptscriptstyle\boldsymbol{\bot}"'] \arrow[d, xshift=0.9ex, "\Uu_{[r+1]}"{name=P3}]\\
\arrow[phantom, from=N, to=Q, , "\scriptscriptstyle\boldsymbol{\dashv}"] 
\arrow[phantom, from=N2, to=Q2, , "\scriptscriptstyle\boldsymbol{\dashv}"] 
\arrow[phantom, from=N3, to=Q3, , "\scriptscriptstyle\boldsymbol{\dashv}"]   
0\text{-}\mathsf{C}\arrow[u, xshift=-0.9ex, "\Ff_{[0]}"{name=F}] 
\arrow[rr,yshift=0.8ex,  "\Tt^r"{name=T2}]&&
\rbc\arrow[u, xshift=-0.9ex, "\Ff_{[r]}"{name=F2}] 
\arrow[rr, yshift=0.8ex,  "\Tt"{name=T2d}]\arrow[ll, yshift=-0.6ex, "(\Tt)^{-r}"{name=U2}]
&&(r+1)\text{-}\mathsf{C} \arrow[ll, yshift=-0.6ex, "\Tt^{-1}"{name=U3}]
\arrow[u, xshift=-0.9ex, "\Ff_{[r+1]}"{name=F3}]   \\
\arrow[phantom, from=P, to=F, , "\scriptscriptstyle\boldsymbol{\dashv}"] 
\arrow[phantom, from=P2, to=F2, , "\scriptscriptstyle\boldsymbol{\dashv}"] 
\arrow[phantom, from=P3, to=F3, , "\scriptscriptstyle\boldsymbol{\dashv}"] 
\arrow[phantom, from=T2, to=U2, , "\scriptscriptstyle\boldsymbol{\sim}"]
\arrow[phantom, from=T2d, to=U3, , "\scriptscriptstyle\boldsymbol{\sim}"] 
    \end{tikzcd}
      \end{center}
      \vspace{-1.5cm}    
\end{cor}

\begin{proof}
The model structures and equivalences of the bottom horizontal line are by Proposition~\ref{P:Qequiv_cxs}. The model structures of
 the middle
horizontal line are by Theorem~\ref{T:trans_model_lwb}. The Quillen equivalences
on the middle horizontal line, as well as the lower vertical Quillen equivalences are by Theorem~\ref{T:eq_lwbr_rbc}.
The model structures of the top line and the upper vertical Quillen equivalences are by
Theorem~\ref{T:modelespse2}.

We check that  the equivalences on the top horizontal line are Quillen equivalences. For
 $r\geq 1$, it is clear that $\Dec(\Ee_{r+1})=\Ee_{r}$. This also holds for $r=0$, using that $H(\rho_1)\colon H((\Nn X)_1)\to H(X_1)$ 
is an isomorphism, by Proposition~\ref{P:Sacyclic}. 
If $f$ is an $(r+1)$-fibration, $\Dec(f)$ is an $r$-fibration using Proposition~\ref{P:Nsurj}.
And the unit of the adjunction
$\eta_X\colon X\to \mathscr{\Dec}\circ \Shift(X)$ is in $\Ee_i$ for all $i\geq 0$,
since it is an isomorphism. So the adjunction $(\Shift,\Dec)$ is a Quillen equivalence.

We now show that $(\LDec, \Shift)$ is a Quillen equivalence. We have adjunctions 
\begin{center}  
    \begin{tikzcd}
\espse_r  
         \arrow[rr, yshift=0.8ex, "\Shift"{name=Wd}] 
&&\espse_{r+1} \arrow[ll, yshift=-0.6ex, "\Dec"{name=Decd}] 
\arrow[rr, yshift=0.8ex, "\LDec"{name=Ld}] 
\arrow[phantom, from=Wd, to=Decd, , "\scriptscriptstyle\boldsymbol{\bot}"]
&&\espse_{r} \arrow[ll, yshift=-0.6ex, "\Shift"{name=S}] 
\arrow[phantom, from=Ld, to=S, , "\scriptscriptstyle\boldsymbol{\bot}"]
\end{tikzcd}
\end{center}
where the left adjunction is a Quillen equivalence. 
The right adjunction is a Quillen adjunction since $\Shift$ preserves fibrations and
weak equivalences.
Since we have $\Dec\circ \Shift\cong 1$ and 
    $\LDec\circ \Shift\cong 1$, the composite adjunction is an adjoint equivalence
    and thus a Quillen equivalence. So the right adjunction is a Quillen equivalence
    by the two-out-of-three property.
\end{proof}

\subsection{Main model category structure on extended spectral sequences}

In this section, for each $r\geq 0$, we prove the existence of a different model structure on extended spectral sequences, with
a stronger notion of weak equivalence, but such that restriction gives the same weak equivalences on spectral sequences.
This structure has the important feature that, in the case $r=0$, spectral sequences sit inside extended spectral sequences as
a homotopically full subcategory. We will explain and exploit this property in Section~\ref{S:infcat}.
The class of  fibrations, $\Fib_r$, is unchanged and so the previous model structure $\espse_r$ is a localization of the one established here.

\begin{defi}
We denote by $\Ee'_r$ the class of maps  $f\colon X\to Y$ of extended spectral sequences 
 satisfying $f_r$ is a quasi-isomorphism and $f_s$ is an isomorphism for $s>r$. Let
 $\Cof'_r$ be  the class of morphisms having the left lifting property with respect to maps in $\Ee'_r\cap\Fib_r$.
\end{defi}

Remark~\ref{R:newequiv} shows that
any map having the left lifting property with respect to $\Fib_r$ is in $\Ee'_r$. We will show directly that there is a model
category with these new equivalences. The proof uses the existence of the model structures $\espse_r$ of
Theorem~\ref{T:modelespse2}.

We write $\Iso_{\leq r}$ for the class of morphisms $f$ in $\espse$ such that $f_i$ is an isomorphism
for $0\leq i\leq r$.

\begin{prop}\label{P:lift}
Morphisms in $\Iso_{\leq r}$ have the left lifting property with respect to those in $\Ee'_r$ and vice versa.
\end{prop}

\begin{proof}
If $u\in \Iso_{\leq r}$ and $p\in \Ee'_r$, there is a lift $h\colon Z\to Y$ in the diagram
\[
	\xymatrix{ X\ar[r]^-{\alpha}\ar[d]_u & Y\ar[d]^{p}\\ Z\ar[r]^{\beta} & W}
\]
given by
\[
h_s=\begin{cases}
	\alpha_s u_s^{-1}&0\leq s\leq r,\\
	p^{-1}_s\beta_s&s\geq r+1,
	\end{cases}
\]
which is compatible with characteristic maps using that $H(p_r)$ is an isomorphism.

If  $u\in \Ee'_r$ and $p\in \Iso_{\leq r}$, there is a lift $h'\colon Z\to Y$ given by
\[
h_s=\begin{cases}
	p^{-1}_s\beta_s&0\leq s\leq r,\\
	\alpha_s u_s^{-1}&s\geq r+1,
	\end{cases}
\]
which is compatible with characteristic maps using that $H(u_r)$ is an isomorphism.
\end{proof}

Note that the lifts are functorial.

\begin{prop}\label{P:cofnew}
Let $f \in\Cof_r$ and let  $g=uf$ where $u\in \Iso_{\leq r}$.
Then $g\in\Cof'_r$.
\end{prop}

\begin{proof}
Consider the left hand diagram with $p\in \Fib_r\cap \Ee'_r$ :
\[
	\xymatrix{ X\ar[r]^-{\alpha}\ar[d]_g & Y\ar[d]^{p}\\ Z\ar[r]_{\beta} & W}
	\qquad \qquad
	\xymatrix{ X\ar[rr]^-{\alpha}\ar[d]_f && Y\ar[d]^{p}\\ Z'\ar[r]_{u}\ar[rru]^{\tilde{h}} &Z\ar[r]_{\beta}\ar[ru]_{h} & W}
\]
The factorization of $g$ as $g=uf$ gives the outer diagram on the right.
 Since $f\in \Cof_r$ and $p\in \Fib_r\cap \Ee'_r\subset \Fib_r\cap \Ee_r$, we have a lift $\tilde{h}\colon Z'\to Y$ and then, since $u\in\Iso_{\leq r}$, we have a lift $h:Z\to Y$ by Proposition~\ref{P:lift}.
\end{proof}

\begin{prop}\label{P:facto3}
In $\espse$ any map $f\colon X\to Y$ with $f\in\Ee_r$ has a functorial factorization $f=\tilde{f}u$ where $u\in \Iso_{\leq r}$ and $\tilde{f} \in \Ee'_r$.
If $f \in \Fib_r$, so is $\tilde{f}$.
\end{prop}

\begin{proof}
Define $Y'$ in $\espse$ by 
\[
	Y'_m=
	\begin{cases}
		X_m&\text{for $m\leq r$,}\\
		Y_m&\text{for $m\geq r+1$}, 
	\end{cases}
\]
with characteristic maps 
\[
	\varphi^{Y'}_m=
	\begin{cases}
		\varphi^X_m&\text{for $m\leq r$,}\\
		(Hf_r)^{-1}\varphi_{r+1}^Y&\text{for $m=r+1$,}\\
		\varphi_m^Y&\text{for $m\geq r+2$.}
	\end{cases}
\]

Then $u\colon X\to Y'$ and  $\tilde{f}\colon Y'\to Y$ given by
\[
	u_m=
	\begin{cases}
		1&\text{for $m\leq r$,}\\
		f_m&\text{for $m\geq r+1$,}
	\end{cases}
	\qquad
	\tilde{f}_m=
		\begin{cases}
		f_m&\text{for $m\leq r$,}\\
 		1&\text{for $m\geq r+1$,} 
		\end{cases}
\]
are compatible with characteristic maps and it is clear that
$u\in  \Iso_{\leq r}$ and $\tilde{f} \in \Ee'_r$ 
Since $\tilde{f}_m=f_m$ for $0\leq m\leq r$, it is clear that if 
$f \in \Fib_r$, so is $\tilde{f}$.
\end{proof}

\begin{prop}\label{P:main-fact}
Any map $f\colon X\to Z$ in $\espse$ has a functorial factorization $f=pi$ where $i\in \Cof'_r$ and $p \in \Fib_{r}\cap \Ee'_r$.
\end{prop}

\begin{proof}
Use the model structure $\espse_{r}$ to factorize $f$ as $f=g\iota$ 
where $\iota\colon X\to Y$ is in $\Cof_{r}$  and $g\colon Y\to Z$ is in $\Fib_{r}\cap \Ee_r$.
Now factorize $g$ as $g=\tilde{g}u$ by Proposition~\ref{P:facto3}, where
 $u\in \Iso_{\leq r}$ and $\tilde{g} \in \Fib_r\cap \Ee'_r$. Then $f=\tilde{g}(u\iota)$
 and $u\iota\in \Cof'_r$ by Proposition~\ref{P:cofnew}.
\end{proof}

\begin{teo}\label{T:modelespse3} 
 There is a model structure on $\espse$ denoted $\espse'_{r}$ with fibrations
$\Fib_r$, the class of maps $f:X\to Y$ such that $f_i$ is bidegreewise surjective for $0\leq i\leq r$, and weak equivalences $\Ee'_r$, the class of maps $f:X\to Y$ such that $f_r$ is a quasi-isomorphism and $f_i$ is an isomorphism for $i>r$.
  \end{teo} 

\begin{proof}
We give a direct proof, checking the axioms as given in~\cite[Definition 2.1.3]{Balchin-handbook}. Firstly it is clear that the weak equivalences, fibrations and cofibrations are closed under composition. Axiom (MC1) is the existence of all small limits and colimits, which holds by Theorem~\ref{T:espse-lims-colims}. (MC2) is the two-out-of-three property for weak equivalences, which is clear.

Axiom (MC3) is that each of the three classes is closed under retracts. This is clear for fibrations and weak equivalences and it holds for
cofibrations as they are defined by the left lifting property.

We consider axiom (MC5) next.
This comprises the two factorization properties: for $f\colon X\to Y$ in $\espse$, we need to
be able to factorize (1)  $f=pi$ where $i\in  \Cof'_r$ and $p\in  \Ee'_r \cap \Fib_r$ and 
(2) $f=p'i'$ where $i'\in \Ee'_r \cap \Cof'_r$ and $p'\in  \Fib_r$.

The first factorization property holds by Proposition~\ref{P:main-fact}. 
Since $\Cof_r\subseteq \Cof'_r$, it follows from Remark~\ref{R:newequiv} that
we have $\Ee_r\cap\Cof_r\subseteq \Ee'_r\cap\Cof'_r$.
So the second factorization property holds because
it holds in the model structure $\espse_r$.

Finally, axiom (MC4) comprises the two lifting properties: given a diagram in $\espse$
\vspace{-0.5cm}
\[
	\xymatrix{ X\ar[r]^-{f}\ar[d]_i & Y\ar[d]^{p}\\ Z\ar[r]^{g} & W} 
\]
there must be a lift (1) when $i\in \Cof'_r$ and $p\in \Ee'_r \cap \Fib_r$ and
(2) when  $i\in \Ee'_r \cap \Cof'_r$ and $p\in  \Fib_r$. We have defined $\Cof'_r$ precisely so that (1) holds.
Lifting property (2) holds by a standard retract argument: factorize $i$ in $\espse_r$ as $i=\rho\iota$ where
$\iota\in\Ee_r\cap\Cof_r$ and $\rho\in \Fib_r$. Then $\iota\in \Ee'_r$ by Remark~\ref{R:newequiv}, so
$\rho\in\Ee'_r$ by the two-out-of-three property. Thus $i$ left lifts against $\rho$ and so $i$ is a retract of $\iota$.
Hence $i\in \Ee_r\cap\Cof_r$ and $i$ left lifts against $p$.
\end{proof}

\subsection{Comparison of model structures}
\label{subsec:compare}

\begin{prop}\label{P:w-dec-espse-new}
The adjunction $\Shift\dashv\Dec$ gives a Quillen equivalence between 
$\espse'_{r}$ and  $\espse'_{r+1}$ for all
$r\geq 0$. And the adjunction $\LDec \dashv \Shift$ gives a Quillen equivalence between 
$\espse'_{r+1}$ and  $\espse'_{r}$ for all
$r\geq 0$.
\end{prop}

\begin{proof}
The proof is essentially the same as that of Corollary~\ref{C:w-dec-espse}, noting for the first
statement that
 $\Dec(\Ee'_{r+1})=\Ee'_{r}$ for all $r\geq 0$ and 
that the unit of the adjunction
$\eta_X\colon X\to \mathscr{\Dec}\circ \Shift(X)$ is in $\Ee'_i$ for all $i\geq 0$, since it is an isomorphism. 
For the second statement, the adjunction is a Quillen adjunction since $\Shift$ preserves fibrations and
weak equivalences. And the same two-out-of-three argument as for Corollary~\ref{C:w-dec-espse} concludes the proof.
\end{proof}

\begin{rmk}\label{R:LDec-preseves-we}
   Note that $\LDec(\Ee'_{r+1})\subset \Ee'_r$ for $r\geq 1$. For $r=0$ this is not true in general. Indeed, let $f:X\to Y$ be in $\Ee'_1$ and consider, using Notation~\ref{N:qr},
   \[\LDec(f)_0\colon\Tt^{-1}\cone_1(X_0)/\Img q_1^X\to \Tt^{-1}\cone_1(Y_0)/\Img q_1^Y.\]  We see that $\LDec(f)_0$ is a quasi-isomorphism if and only if the induced map $\Ker q_1^X\to\Ker q_1^Y$ is a quasi-isomorphism. We know this is true for $X$ and $Y$ cofibrant in $\espse'_{1}$, because $\LDec$ is a left Quillen functor. It is also true if $X$ and $Y$ are spectral sequences because then $\Ker q_1^X=\Ker q_1^Y=0$.
   \end{rmk}

The following result explains the relationship between the new model structures and those obtained in the previous section.

\begin{prop}\label{P:localization}
The identity functor $\espse'_r\to \espse_r$ is a right Bousfield localization which is not a Quillen equivalence.
\end{prop}

\begin{proof}
The first part is immediate since we have established the two model category structures on $\espse$ with the same
fibrations and $\Ee'_{r}\subset \Ee_r$.

Let $X$ be a cofibrant object of  $\espse_r$ and $Y$ any object in $\espse$.
If this were a Quillen equivalence, we would have that $f\colon X\to Y$ is in $\Ee_r$ if and only if $f\colon X\to Y$ is in $\Ee'_{r}$. Since $0$ is cofibrant 
this implies that any extended spectral sequence $X$ such that $H(X_r)=0$ satisfies $X_i=0$ for all $i>r$ which is certainly not true.
\end{proof}

So all our model structures are Quillen equivalent to one of either $\espse_0$ or $\espse'_0$. 
Since $\espse_0$ is Quillen equivalent to the projective model structure on $0$-bigraded complexes, Proposition~\ref{P:localization} can be 
phrased as saying that $\espse'_0$ is a right \emph{delocalization}
of (a model category Quillen equivalent to) the projective model structure on $0$-bigraded complexes.

\subsection{Cofibrations in $\espse_0'$}\label{subsec:cofibs}
In this section we characterize the cofibrations of the model category $\espse_0'$.

\begin{lem}\label{L:lift0page}
Let $f\colon X \to Z$ be a map in $\espse$.
Any commutative diagram as on the left
\begin{equation}
\begin{gathered}
	\xymatrix{ X_0\ar[r]^-{\alpha}\ar[d]_{f_0} & Y\ar[d]^{p}\\ Z_0\ar[r]^{\beta} & W}
\qquad \qquad
\xymatrix{ X\ar[r]^-{\overline{\alpha}}\ar[d]_{f} & \overline{Y}\ar[d]^{\overline{p}}\\ Z\ar[r]^{\overline{\beta}} & \overline{W}}
\end{gathered}
\tag{$\ast$}
\end{equation}
in $\widehat{\langle\delta_0\rangle}$ is obtained as the $0$-page of a commutative diagram in $\espse$
as on the right. Furthermore, if $p$ is a surjective quasi-isomorphism, then $\overline{p}\in \Fib_0\cap\Ee_0$.
\end{lem}

\begin{proof}
A diagram with the required properties is obtained by setting
\[
	\overline{Y}_m=
	\begin{cases}
		Y&\text{for $m=0$,}\\
		X_m&\text{for $m\geq 1$,}
	\end{cases}
	\qquad
	\varphi^{\overline{Y}}_m=
		\begin{cases}
		(H\alpha)\varphi_1^X&\text{for $m=1$,}\\
 		\varphi_m^X&\text{for $m\geq 2$,} 
		\end{cases}
\]
\[
	\overline{W}_m=
	\begin{cases}
		W&\text{for $m=0$,}\\
		Z_m&\text{for $m\geq 1$,}
	\end{cases}
	\qquad
	\varphi^{\overline{W}}_m=
		\begin{cases}
		(H\beta)\varphi_1^Z&\text{for $m=1$,}\\
 		\varphi_m^Z&\text{for $m\geq 2$,} 
		\end{cases}
\]
\[
	\overline{\alpha}_m=
	\begin{cases}
		\alpha&\text{for $m=0$,}\\
		1&\text{for $m\geq1$,}
	\end{cases}
	\qquad
	\overline{\beta}_m=
	\begin{cases}
		\beta&\text{for $m=0$,}\\
		1&\text{for $m\geq1$,}
	\end{cases}
	\qquad
	\overline{p}_m=
	\begin{cases}
		p&\text{for $m=0$,}\\
		f_m&\text{for $m\geq1$.}
	\end{cases}
\]
It is routine to check that $\overline{\alpha}, \overline{\beta}, \overline{p}$ are maps of extended spectral sequences.
It is clear that the diagram commutes and that taking the $0$-page recovers the original diagram. The final statement is immediate
since $\overline{p}_0=p$. 
\end{proof}

\begin{prop}\label{P:new0cof}
The cofibrations $\Cof'_0$ in $\espse'_0$ are those $f\colon X\to Y$ such that
 $f_0$ is a cofibration of $0$-bigraded complexes in the model category structure on $0\text{-}\mathsf{bC}$ given in Proposition~\ref{P:modcatrbc}.
 That is, $f\in\Cof'_0$ if and only if $f_0$ has the left lifting property with respect to surjective quasi-isomorphisms of $0$-bigraded complexes.
\end{prop}

\begin{proof}
Suppose that $f$ is such that $f_0$ has the left lifting property with respect to surjective quasi-isomorphisms of $0$-bigraded complexes.
Consider the diagram in $\espse$
\[
	\xymatrix{ X\ar[r]^-{\alpha}\ar[d]_f & Z\ar[d]^{p}\\ Y\ar[r]^{\beta} & W}
\]
where $p\in \Fib_0\cap \Ee'_0$.
There is a lift $h_0\colon Y_0\to Z_0$ on the $0$-page by assumption and then
\[
h_s=\begin{cases}
	h_0&s=0,\\
	p_s^{-1}\beta_s&s\geq 1,
\end{cases}
\]
gives a lift, using that $H(p_0)$ is an isomorphism for compatibility with characteristic maps.

Conversely, let $f\in\Cof'_0$ and let $p$ be a surjective quasi-isomorphism of $0$-bigraded complexes. We need to show there
is a lift in the diagram appearing on the left of $(\ast)$ in the statement of Lemma~\ref{L:lift0page}
when $p$ is a surjective quasi-isomorphism.
By the lemma we have the diagram on the right of $(\ast)$
and $\overline{p}\in \Fib_0\cap\Ee_0$.

Now apply Lemma~\ref{P:facto3} to factorize $\overline{p}$ as $\overline{p}=\tilde{p}u$ where $u\colon\overline{Z}\to Z$  is 
in $\Iso_{\leq 0}$ and 
$\tilde{p}\colon Z\to\overline{W}$ is in $\Fib_0\cap \Ee'_0$. Since $f\in\Cof'_0$ there is a lift $h\colon Y\to Z$ 
\[
	\xymatrix{ X\ar[r]^-{\overline{\alpha}}\ar[d]_f & \overline{Z}\ar[r]^{u}&Z\ar[d]^{\tilde{p}}\\ 
	Y\ar[rr]^{\overline{\beta}}\ar[rru]^{h} &&\overline{W}  }
\]
and then $u_0^{-1}h_0$ provides the
required lift.
\end{proof}

\section{Infinity-categorical interpretation}\label{S:infcat}
\setcounter{subsection}{1} 

In this section we will interpret our results in the context of relative categories in the setting of Barwick-Kan~\cite{BK}, as well as
explaining the relationship to our previous work~\cite{LW24}.

In~\cite[Theorem 5.3.1]{LW24} we established the existence of structures on the category of spectral sequences $\spse$ 
closely related to the model category structures on $\espse$ presented here. The relevant structure has a class of fibrations and a class of weak equivalences and is called an \emph{almost Brown category with functorial factorization}. Given $r\geq 0$, there is such a structure,
denoted $\spse_r$, with fibrations the class $\Fib_r$ of morphisms of spectral sequences $f$ such that $f_i$ is surjective for $0\leq i\leq r$ and weak equivalences $f$ such that $H(f_r)$ is an 
isomorphism. 

\begin{defi}
A \emph{relative category} is a pair $(\mathcal{C}, \mathcal{W})$ where $\mathcal{C}$ is a category and
$\mathcal{W}$ is a subcategory containing all the objects of $\mathcal{C}$; the maps in $\mathcal{W}$ are called
\emph{weak equivalences}.
We write  $\RelCat$  for the category of relative categories endowed with the model category structure of Barwick-Kan~\cite{BK}.
\end{defi}

We study the relative category $(\espse, \Ee_r')$ and its restriction 
$(\spse, \Ee_r)$.

We have seen that for $f\colon X\to Y$ a morphism in $\espse$ we have $f\in\Ee'_{r+1}$ if and only if $\Dec f\in\Ee'_r$. And it is easy to see that
$f\in \Ee'_r$ if and only if $\Shift f \in \Ee'_{r+1}$.

We will use ideas and results of Pascaleff~\cite{Pas23} and the main theorem of Meier~\cite{Meier16}.

\begin{defi}{\cite[Definition 2.1]{Pas23}}
A \emph{Dwyer-Kan adjunction} between relative categories is an adjunction 
$(L,R)\colon (\Ccat_1,\Ww_1)\to (\Ccat_2,\Ww_2)$ such that $L$ and $R$ preserve the class of weak equivalences 
and such that the unit of the adjunction lies in $\Ww_1$ and the counit of the adjunction in $\Ww_2$.
\end{defi}

As noted by Pascaleff, a result of Bousfield-Kan~\cite[Corollary 3.6]{DK80} shows that if two 
relative categories are related by a Dwyer-Kan adjunction, then they are Barwick-Kan equivalent in $\RelCat$.

Consider the following commutative diagram of relative categories.

 \begin{center}  
    \begin{tikzcd}[row sep=8ex, column sep=8ex]
(\espse, \Ee'_0) \arrow[r, yshift=0.8ex, "\Shift^r"{name=W}]
&(\espse, \Ee'_r) \arrow[l, yshift=-0.6ex, "\Dec^r"{name=Dec}]  
\arrow[r, yshift=0.8ex, "\Shift"{name=Wd}]&(\espse, \Ee'_{r+1}) 
\arrow[l, yshift=-0.6ex, "\Dec"{name=Decd}]  \\
\arrow[phantom, from=W, to=Dec, , "\scriptscriptstyle\boldsymbol{\bot}"] 
\arrow[phantom, from=Wd, to=Decd, , "\scriptscriptstyle\boldsymbol{\bot}"]
(\spse, \Ee_{0}) \arrow[u, hook] 
\arrow[r, "\Shift^r" description]
&(\spse, \Ee_r)\arrow[u, hook]  
\arrow[l, yshift=2.5ex, "\LDec^r"'{name=LDec2d}, "\scriptscriptstyle\boldsymbol{\bot}"]
\arrow[l, yshift=-2.5ex, "\Dec^r"{name=Dec2}, "\scriptscriptstyle\boldsymbol{\bot}"'] 
\arrow[r, "\Shift" description]
&(\spse, \Ee_{r+1})\arrow[u, hook].   
\arrow[l, yshift=2.5ex, "\LDec"'{name=LDec2d}, "\scriptscriptstyle\boldsymbol{\bot}"]
\arrow[l, yshift=-2.5ex, "\Dec"{name=Dec2d}, "\scriptscriptstyle\boldsymbol{\bot}"']
    \end{tikzcd}
      \end{center}   
      
      \begin{teo}
\label{T:fib-recat}
 The adjunctions in the diagram above are Dwyer-Kan adjunctions. In particular the rows are Barwick-Kan equivalences in $\RelCat$. In addition the relative categories $(\spse, \Ee_0)$  and  $(\espse, \Ee'_r)$ for $r\geq 0$ are fibrant relative categories in the Barwick-Kan model structure on $\RelCat$.
      \end{teo}
      
      \begin{proof} We consider the $(\Shift, \Dec)$ adjunction.
      The unit  is an isomorphism.  The counit  $\epsilon\colon \Shift\circ\Dec\to \Id_{\espse}$ is the identity at page $i\geq 2$ and a quasi-isomorphism at page 1, using Proposition~\ref{P:Sacyclic}. Hence the counit lives in $\Ee'_{r}$ for every $r\geq 1$.   
     The category $\espse'_r$ is a model category and so 
     $(\espse, \Ee'_r)$ is
     fibrant by~\cite[Theorem 4.13]{Meier16}.
Furthermore if $f\colon X\to Y$ is in $\Ee'_{0}$ where one of $X$ or $Y$ is a spectral sequence, then so is the other. Hence $(\spse_0, \Ee_0)$ is a relative category which is  homotopically full in $(\espse, \Ee'_0)$, in the sense of Meier and so it is fibrant by~\cite[Theorem 4.13]{Meier16}.

Now we consider the $(\LDec, \Shift)$ adjunction on $\spse$. We have seen that
$\Shift$ preserves weak equivalences and by Remark~\ref{R:LDec-preseves-we}, 
$\LDec\colon (\spse, \Ee_{r+1})\to (\spse, \Ee_r)$ preserves weak equivalences for all $r$. The counit of the 
adjunction is an isomorphism. The unit of the adjunction is the identity at pages $i\geq 2$ and, by the same reasoning as in
Remark~\ref{R:LDec-preseves-we}, it is
a quasi-isomorphism at page $1$ on $\spse$ (but not in general on $\espse$).
          \end{proof}

\begin{cor}
    For each $r\geq 1$, the relative category
$(\spse, \Ee_r)$ has $(\spse, \Ee_0)$ as a fibrant replacement
in the Barwick-Kan model structure on $\RelCat$, via the functor $\Dec^r$, or the functor $\LDec^r$. \qed
\end{cor}

\newpage

\appendix
\section{Representable linear witness books}\label{App:liftingprops}

\subsection{Description of the representable linear witness book 
\texorpdfstring{$\Yy(r,p,n)$}{Y}}
\begin{example}\label{E:Yrpn} We treat the combinatorics of the representable object $\Yy(r,p,n)$ of $\lwb$. This
 is the free linear witness book generated by an element $\alpha_r^{p,n}$, corresponding to  the identity map in $\Hom_{\D}((r,p,n),(r,p,n))$.
 
We describe it in  terms of Proposition~\ref{P:hats}. For each $i$, the $i$-bigraded complex $(\Yy(r,p,n))_i$ has two bigraded components free of rank one, and the other components are $0$. 
More precisely we have
\begin{align*}
\Yy(r,p,n)_i&=\kk\;  (w_r)^{r-i}d_r\alpha_r^{p,n}\oplus\kk\; (w_r)^{r-i}\alpha_r^{p,n}\quad\text{for $0\le i <r$,}\\
\Yy(r,p,n)_r&=\kk\; d_r\alpha_r^{p,n}\oplus\kk\; \alpha_r^{p,n},\\
\Yy(r,p,n)_{r+k}&=\kk\; d_{r+k}(s_r)^k\alpha_r^{p,n}\oplus\kk\; (s_r)^k\alpha_r^{p,n}\qquad\text{for $k\geq 0$}. 
\end{align*}
And the actions of $d,s,w$ are clearly identified knowing the relations of Proposition \ref{P:hats} (compare also with Lemma \ref{L:morinS}).
Here we denote by $(w_r)^j$ the composite $w_{r-j+1}\ldots w_{r+1}w_r$ and the same notation is used for the 
iterated composite of the $s$s. 

In addition, any $\psi\colon (r,p,n)\rightarrow (r',p',n')$ induces a morphism of linear witness books $\Yy(r,p,n)\rightarrow \Yy(r',p',n')$ sending the identity morphism to the morphism $\psi$ in $\Yy(r',p',n')$. More concretely we have
\[\xymatrix@R=0.5pt@C=2pc{\Yy(r,p,n)\ar[r]^-{\delta_r} & \Yy(r,p+r,n+r-1)& \Yy(r,p,n)\ar[r]^-{\omega} & \Yy(r+1,p,n)\\
 \varphi \alpha_r^{p,n}\quad \ar@{|->}[r] & \quad\varphi d_r \alpha_{r}^{p+r,n+r-1}&\varphi \alpha_r^{p,n}\quad \ar@{|->}[r] & \quad\varphi w_{r+1} \alpha_{r+1}^{p,n}
}\]
\[\xymatrix@R=0.5pt@C=2pc{\Yy(r,p,n)\ar[r]^-{\sigma} & \Yy(r-1,p-1,n-1)\\
 \varphi \alpha_r^{p,n}\quad \ar@{|->}[r] & \quad\varphi s_{r-1} \alpha_{r-1}^{p-1,n-1}.
 }\]
 \smallskip
 
Of course, $ \Hom_{\lwb}  (\Yy(r,p,n),A)\cong A(r,p,n)$.
\end{example}

\subsection{Characterizations of $\lwbe$ and $\lwbs$ via lifting properties}

\begin{notation}\label{N:rep_ker} Let $\Zz(r,p,n)$ (resp. $\Sf(r,p,n))\in\lwb$ be the representative of  $\Ker d_r$ (resp. $\Ker s_r$). These objects are defined as the following  cokernels.
\begin{align*}
\Zz(r,p,n)& =\mathrm{coker}(\delta_r\colon \Yy(r,p-r,n+1-r)\rightarrow \Yy(r,p,n)),\\
\Sf(r,p,n)& =\mathrm{coker}(\sigma_r\colon \Yy(r+1,p+1,n+1)\rightarrow \Yy(r,p,n)).
\end{align*}

Since any linear witness book satisfies $\Ker s_r\subseteq \Ker d_r$, there is a (surjective) morphism 
$\tau_r\colon \Zz(r,p,n)\to \Sf(r,p,n)$.  Let $M^\tau$ denote the set of morphisms of $\lwb$ given by 
\[
M^\tau\colon =\left\{\tau_r\colon  \Zz(r,p,n)\rightarrow \Sf(r,p,n) \right\}_{r\geq 0,\, p,n\in\ZZ}.
\]

We denote by $\kk_{\leq r}(p,n)$ the extended spectral sequence having $\kk$ concentrated in bidegree $(p,n)$ at pages $i$ such that $0\leq i\leq r$ and $0$ elsewhere, with characteristic maps being the identity or zero.
\end{notation}

\begin{lem}\label{L:Z_notin_lwbe} The linear witness book $\Sf(r,p,n)$ lies in $\lwbe$, but the linear witness book $\Zz(r,p,n)$ does not. We have the following isomorphisms in $\lwbe$: 
\[
	\Sf(r,p,n)\cong\Nn\Qq\Sf(r,p,n)=\Nn\Qq\Zz(r,p,n)\cong\Nn(\kk_{\leq r}(p,n)).
\]
\end{lem}

\begin{proof} We first compute $\Sf(r,p,n)$ and $\Zz(r,p,n)$ in $\lwb$, thanks to Example \ref{E:Yrpn}.
On the one hand we see that at page $i\leq r$, $\Sf(r,p,n)_i=\Zz(r,p,n)_i=\kk w^{r-i}{\alpha}_r^{p,n}$, with $d_i=0$, $w_i$ is the 
identity and  $s_i=0$, for $i<r$.
At page $i>r$ we have $\Sf(r,p,n)_i=0$. This implies $\Ker d_i=\Ker s_i$ in $\Sf(r,p,n)_i$, for all $i\geq 0$, so that $\Sf(r,p,n)\in\lwbe$.
At page $i>r$ we have $\Zz(r,p,n)_i=\kk d_i s^{i-r}{\alpha}_r^{p,n}\oplus \kk s^{i-r}{\alpha}_r^{p,n}$. In particular $\Ker s_r=0$ while $\Ker d_r=\kk\alpha_r^{p,n}$.
In addition we have $\Qq\Sf(r,p,n)=\Qq\Zz(r,p,n)=\kk_{\leq r}(p,n)$.
\end{proof}

The following characterization of $\lwbe$ via right lifting properties is an immediate consequence of the descriptions.

\begin{prop}\label{P: charespse_rlp}
\begin{enumerate}
\item A map $f\colon K\to L$ has the right lifting property with respect to all morphisms in $M^\tau$ if and only if
whenever $a\in \Ker d_r$ and $f(a)\in \Ker s_r$ then $a\in \Ker s_r$.
\item A linear witness book $L$ is in $\lwbe$ if and only if $L\to 0$ has the right lifting property with respect to
morphisms in $M^\tau$. \qed
\end{enumerate}
\end{prop}

\begin{rmk}\label{local-lwbe}
Equivalently, we can characterize $\lwbe$ as the subcategory of $M^\tau$-local objects in $\lwb$, in the ordinary categorical sense, see for example
~\cite[p143]{Riehl_context}. That is, an $M^\tau$-local object $L$ is one such that for all morphisms $\tau_r\colon \Zz\to\Sf$ in $M^\tau$, the induced map $\lwb(\Sf, L)\to
\lwb(\Zz, L)$ is a bijection. Surjectivity is by Proposition~\ref{P: charespse_rlp} and injectivity holds because of surjectivity of the maps
$\tau_r$ in $M^\tau$.
\end{rmk}

The characterization of objects of $\lwbs$ in Proposition~\ref{P: charspseq} can also be described by means of right lifting properties. Namely, for $r\geq 1$,
let 
\[\Ww(r,p,n)=\Coker(\omega_{r}\colon \Yy(r-1,p,n)\rightarrow \Yy(r,p,n)),\]
representing $\Ker w_r$.
Note that $\Ww(r,p,n)$ is in $\lwbe$ and that the extended spectral sequence $\Qq\Ww(r,p,n)$ takes the following form:
\[\Qq\Ww(r,p,n)_i=\begin{cases} \kk w^{r-i}d_r\alpha_r^{p,n}&\text{  if }\ 0\leq i<r,\\
\kk d_r\alpha_r^{p,n}\oplus\kk \alpha_r^{p,n}&\text{  if }\ i=r,\\
0&\text{  if }\ i>r.\end{cases}\]

The relations $\omega_{r+1}\sigma_r=0$ and $\sigma_r\omega_{r+1}=0$ imply that there are well defined sets of morphisms of $\lwb$: 
\begin{align*}
M^{\sigma,\omega}\colon &=\left\{\sigma_r\colon  \Ww(r+1,p+1,n+1)\rightarrow \Yy(r,p,n) \right\}_{r\geq 0,\, p,n\in\ZZ},\\
M^{\omega,\sigma}\colon &=\left\{\omega_{r+1}\colon   \Sf(r,p,n)\rightarrow \Yy(r+1,p,n) \right\}_{r\geq 0,\, p,n\in\ZZ}.
\end{align*}

The following proposition is a reformulation of Proposition~\ref{P: charspseq}.
\begin{prop}\label{P:charspse_rlp}
Let $L$ be an object of $\lwbe$. 
\begin{enumerate}
\item The characteristic maps $\varphi_i$ of $\Qq L$ are surjective for all $i$ if and only if $L\to 0$ has the right lifting property with respect to $M^{\omega,\sigma}$.
\item The characteristic maps $\varphi_i$ of $\Qq L$ are injective for all $i$ if and only if $L\to 0$ has the right lifting property with respect to $M^{\sigma,\omega}$. \qed
\end{enumerate}
\end{prop}

We collect information about key objects in $\lwb$ in the following table. 

\begin{table}[h]
\caption{Key objects in $\lwb$}
\begin{tabular}{ccll}
object&represents&homology\\
\hline
$\Yy(r,p,n)$&element&$\kk^{p,n}\oplus\kk^{p-r,n+1-r}$&on pages $0\le i <r$\\
&&$0$&on pages $i\geq r$\\
\hline
$\Zz(r,p,n)$&$\Ker d_r$
&$\kk^{p,n}$&on pages $0\le i \le r$\\
&&$0$&on pages $i> r$\\
\hline
$\Sf(r,p,n)$&$\Ker s_r$
&$\kk^{p,n}$&on pages $0\le i \le r$\\
&&$0$&on pages $i> r$\\
\hline
$\Ww(r,p,n)$&$\Ker \omega_r$&$\kk^{p-r,n+1-r}$&on pages $0\le i <r$\\
&&$0$&on pages $i\geq r$\\
\hline
\end{tabular}
\end{table}

\section{On the functor $\LDec$}
\label{App:adjoints}

\setcounter{subsection}{1} 

In this section we provide results related to the 
functor \(\LDec\), left adjoint to the shift functor.
The main results are Proposition~\ref{P:ldecr-shiftr-espse}, giving an explicit
description of \(\LDec\) and its iterates, and Proposition~\ref{prop: ldec-fcx},
which shows compatibility with Deligne's left adjoint
$\Dec^*$ to shift
on filtered complexes.

We begin with a description of the functor $(\Ww_{\geq r})_!$, left adjoint to the functor $\Ff_{\geq r}\colon \lwb_{\geq r}\to\lwb.$ This description is obtained by
applying the (co)end formulas of Proposition~\ref{P:Kelly}.

\begin{prop}\label{P:wexclamation}
    The left adjoint functor $(\Ww_{\geq r})_!:\lwb\to\tr_{\geq r}\lwb$ to the functor $\Ff_{\geq r}$ has the following form. For $(L,d_i,w_i,s_i)$ a linear witness book, and for $i\geq r$,
    \[(\Ww_{\geq r})_!(L)_i^{p,n}=(L_0^{p+r-i,n+r-i}\oplus L_0^{p+r;n+r-1}\oplus L_i^{p,n})/\sim\]
    where $\sim$ is the equivalence relation generated by, for all $x\in L_r$,
    \begin{align*}
        (0,0,s^{i-r}x)&\sim(w^rx,0,0),\\
        (0,0,d_i s^{i-r} x)&\sim (0,w^r x,0).
    \end{align*}
    The structure maps $d_i^!,s_i^!,w_{i+1}^!$ are defined for $i\geq r$ by
 \begin{align*}
         d_i^!&=\begin{pmatrix}0&0&0\\1& 0&0\\0&0&d_i\end{pmatrix}, & w_{i+1}^!&=\begin{pmatrix}0&0&0\\0& 1&0\\ 0&0&w_{i+1}\end{pmatrix}, & s_i^!&=\begin{pmatrix}1&0&0\\0& 0&0\\0&0&s_i\end{pmatrix}. 
 \end{align*}
 The counit of the adjunction is an isomorphism at each page $i\geq r$ and the unit of the adjunction $\eta_L:L\to \Ff_{\geq r}(\Ww_{\geq r})_! L$ has the following description.

 \begin{itemize}
     \item For $x\in L_i$ and $0\leq i<r$ 
     \[\eta_L(x)=\overline{(w ^ix,0,0)}\in (\Ww_{\geq r})_!(L)_r.
     \]
     \item For $x\in L_i$ and $i\geq r$
     \[\eta_L(x)=\overline{(0,0,x)} \in (\Ww_{\geq r})_!(L)_i.
     \vspace{-0.75cm}
     \] \qed
 \end{itemize}
    \end{prop}

    Using the characterization of $\lwbe$ in Proposition~\ref{prop:char-espse} and that of $\lwbs$ in Proposition~\ref{P: charspseq} it is immediate to prove the next statement.

\begin{prop}\label{P:wwrestrict}
    The functor $(\Ww_{\geq r})_!$ restricts to a functor $\lwbe\to \tr_{\geq r}\lwbe$ and to a functor $\lwbs\to  \tr_{\geq r}\lwbs$.\qed
\end{prop}

In order to describe $\LDec^r$ explicitly, we introduce some notation.

  \begin{notation}\label{N:qr} Let $X$ be in $\espse$. Denote by $q_r\colon (\Nn X)_r^{p+r,n+r-1}\to \cone_r(X_0)^{p,n}$ the map defined by
  \[
  q_r(x_0,\ldots,x_r;y_0,\ldots,y_r)=(-y_0,x_0).
  \]
  For an $r$-bigraded complex $A$ we denote by $\Sigma A$ the $r$-bigraded complex defined as $(\Sigma A)^{p,n}=(\Sigma A)^{p+r,n+r-1}$ with differential $d^{\Sigma A}=-d^A$.
  The map $q_r$ induces  a short exact sequence in $\rbc$
  \[0\to \Sigma(\Nn X)_r/\Ker q_r\to \cone_r(X_0)\to  \cone_r(X_0)/\Img q_r\to 0\]
  which yields a long exact sequence in homology and thus an isomorphism
  \[\partial_r\colon H(\cone_r(X_0)/\Img q_r)\to H((\Nn X)_r/\Ker q_r).\]
  \end{notation}

Now we give the description of $\LDec^r$.

\begin{prop}\label{P:ldecr-shiftr-espse}
    For $X$ in $\espse$, we have in $\tr_{\geq r}\espse$
     \[(\Tt^r\LDec^r X)_i=\begin{cases}
 \cone_r(X_0)/\Img q_r & \text{ for }  i=r,\\
(X_i,d_i) & \text{ for } r< i.\\
    \end{cases}\]
   The characteristic maps $\varphi_i$ coincide with $\varphi_i^X$ for $r+1<i$ and $\varphi_{r+1}$ 
   is the following composite
 \[
 \begin{tikzcd}[sep=small]
X_{r+1}\arrow[r,"\varphi^X_{r+1}"]& H(X_r)\arrow[r,"H(\rho_r)^{-1}"] & H((\Nn X)_r)\arrow[r] & H((\Nn X)_r/\Ker q_r)
\arrow[r,"\partial_r^{-1}"] & H(\cone_r(X_0)/\Img q_r),
   \end{tikzcd}
   \]
where the unlabelled arrow is induced by the projection.

The counit of the adjunction
$(\LDec^r,\Shift^r)$ is an isomorphism.
    \end{prop}

\begin{proof} Let $X\in\espse$. By Theorem \ref{T:shiftdectheorem} we have
\[\Tt^r\LDec^r(X)\cong \Qq_{\geq r}(\Ww_{\geq r})_!\Nn X.\]
We first observe that for $L\in \lwb$ and $i\geq r$, the $i$-bigraded complex $(\Ww_{\geq r})_!(L)_i$ is a quotient of the $i$-bigraded complex 
$M_i=\cone_i(L_0)^{*+r-i,*+r-i}\oplus L_i$. \\
Assume  $i>r$. A short computation shows that the inclusion $L_i\to M_i$ induces 
 an isomorphism of $i$-bigraded complexes
\[\Qq(L)_i\to \Qq_{\geq r}((\Ww_{\geq r})_!(L))_i.\]
In particular if $L=\Nn X$, then by Theorem \ref{T:adj_QN} we have an isomorphism of $i$-bigraded complexes
\[\Qq_{\geq r}((\Ww_{\geq r})_!(\Nn X))_i\cong (\Qq\Nn X)_i\cong X_i,\]
and one can check that the  characteristic maps $\varphi_i$ coincide with $\varphi_i^X$.\\
Assume $i=r$. We have
\[\Qq_{\geq r}((\Ww_{\geq r})_!(\Nn X))_r\cong((\Ww_{\geq r})_!(\Nn X))_r\cong \cone_r(X_0)\oplus (\Nn X)_r/\sim\]
where, by Proposition \ref{P:wexclamation}, $\sim$ is the equivalence relation generated by
\[(0,0,x)\sim (w^rx,0,0) \text{ and } (0,0,d_rx)\sim (0,w^rx,0),\]
for $x\in(\Nn X)_r$. Let us write $x=(u_0,\ldots,u_r;v_0,\ldots,v_r)$. Recall from Proposition \ref{P:nerve_descr} that $w^rx=u_0$ and $d_rx=(v_0,\ldots,v_r;0,\ldots,0)$.
A quick check shows that the obvious map $\cone_r(X_0)\to \Qq_{\geq r}((\Ww_{\geq r})_!(\Nn X))_r$ is surjective with kernel $\Img q_r$. The expression of $\varphi_{r+1}$ is obtained through these isomorphisms. 

Let $X\in\tr_{\geq r}\espse$ and
let us consider the counit of the adjunction 
\[
(\Tt^r\LDec^r,\Shift^r\Tt^{-r}):\espse\to\tr_{\geq r}\espse.
\]  It is clear that this counit is an isomorphism at pages $i>r$.
At page $r$, it takes the following form:  it sends $\overline{(a,b)}\in \cone_r(X_r)/\Img q_r$ to $a+d_rb\in X_r$ and it is an isomorphism.
 \end{proof}

 The next result establishes the compatibility of $\LDec$ with
 Deligne's left adjoint $\Dec^*$ to shift on filtered complexes.

\begin{prop}
\label{prop: ldec-fcx}
Let $C$ be a filtered complex.
There is a natural isomorphism of spectral sequences $E(\Dec^*(C))\to \LDec(E(C))$.
\end{prop}

\begin{proof}
From~\cite[1.3.5]{DeHII},
we have a canonical morphism of spectral sequences 
\[v\colon \Tt^{-1}\circ \tr_{\geq 1}(E(C)) \to E(\Dec^*(C))\] 
which is an injection on the $0$-page and an isomorphism from the $1$-page onwards.

In fact $v_0\colon  \Tt^{-1}\circ \tr_{\geq 1}(E(C))_0 \to
E(\Dec^*(C))_0$ factors via  $\LDec(E(C))_0$
and we use this to produce the required isomorphism. 
Recall that 
\[
E(\Dec^*(C))_0^{p,p+n}=
\Tt^{-1}(B_1^{p,p+n}(C)/B_1^{p-1,p+n-1}(C))
\]
with notation as in~\cite[Section 2]{CELW20}. Using the
description of $\LDec$ from Proposition~\ref{P:ldecr-shiftr-espse},
it is routine to check that there is a well-defined map
\[\tilde{v}_0\colon  \LDec(E(C))_0=
\Tt^{-1}(\cone_1(E(C)_0)/\Img q_1)
\to E(\Dec^*(C))_0\]
given by
\[
\overline{(a,b)}\mapsto [a+db].
\]
Furthermore, $v_0$ is the composite of the map sending $a$ to $\overline{(a,0)}$ and $\tilde{v}_0$. It is straightforward to check that
$\tilde{v}=(\tilde{v}_0, v_1, v_2, \dots)\colon E(\Dec^*(C))\to \LDec(E(C))$ is an isomorphism of spectral sequences. 
\end{proof}

\medskip

For completeness, we end this section with
the description of the right adjoint functor $\Rr_{[r,t]}$ to the forgetful functor $\Uu_{[r,t]}\colon \lwb\to\tr_{[r,t]}\lwb$. 

\begin{prop}\label{P:RU-LWB} Let $(L,d_i,w_i,s_i)$ be an object of $\tr_{[r,t]}\lwb$.
We have
 \[ \Rr_{[r,t]}(L)_i^{p,n}=\begin{cases}
 L_r^{p+r-i,n+r-i} & \text{ for } 0\leq i<r, \\
 L_i^{p,n}& \text{ for } r\leq i\leq t,\\
 \cone_i(Z_t(L_t))^{p-i,n-i+1}& \text{ for } i>t, \end{cases}\]

with the maps $d_i^\Rr,w_i^\Rr,s_i^\Rr$ defined as
 \begin{align*}
     d_i^\Rr&=0& w_{i+1}^\Rr&=0 & s_i^\Rr&=1 & \text{ for } 0\leq i<r,\\
     d_i^\Rr&=d_i &  w_{i+1}^\Rr&=w_{i+1} & s_i^\Rr&=s_i & \text{ for } r\leq i< t, \\
     d_t^\Rr&=d_t^{L} &  w_{t+1}^\Rr&=\begin{pmatrix}0&1\end{pmatrix} & 
     s_t^\Rr&=\begin{pmatrix}d_t\\0\end{pmatrix} & \text{ for } i=t, \\
     d_i^\Rr&=\begin{pmatrix}0&0\\1& 0\end{pmatrix} & w_{i+1}^\Rr&=\begin{pmatrix}0&0\\0& 1\end{pmatrix} & s_i^\Rr&=\begin{pmatrix}1&0\\0& 0\end{pmatrix} & \text{ for } i>t.
 \end{align*}

 The counit of the adjunction $(\Uu_{[r,t]},\Rr_{[r,t]})$ is the identity. 
 The unit of the adjunction 
$L\rightarrow \Rr_{[r,t]}\Uu_{[r,t]} L$ takes the following form.
\begin{itemize}
\item For $0\leq i<r$ it corresponds to $s^{r-i}\colon L_i^{p,n}\rightarrow L_r^{p+r-i,n+r-i}$.
\item It is the identity at the $i$-page, for $r\leq i\leq t$.
\item For $i>t$ the unit corresponds to  the map  
\[\begin{pmatrix} 0\\w^{i-t}\end{pmatrix}\colon L_i^{p,n}\to Z_t(L_t)^{p-i,n+1-i}\oplus Z_t(L_t)^{p,n}.\vspace{-0.75cm}\]
\qed
\end{itemize}
\end{prop}

\medskip

\begin{rmk}
    If $r>0$ then $\Rr_{\geq r}$ does not restrict to a functor from $\tr_{\geq r}\lwbe$ to $\lwbe$, since $\Ker d_i\not=\Ker s_i$ for $i<r$.
\end{rmk}

\section{Non-existence of certain model categories}
\label{App:nonmodels}
\setcounter{subsection}{1} 

In this section, we prove that there is no model category structure on either $\lwb$ or on $\lwbe$ such that $\lwbs$ is the subcategory of fibrant objects, where  the class of weak equivalences between fibrant objects   is the class $\Ee_r$ (or equivalently $\Ee'_r$). 

We characterize morphisms $p\colon L\to M$ in $\lwb$ with $L,M$ in $\lwbs$, such that $\Ker p$ lives in $\lwbs$. Recall that $\Ker p$ is computed in $\lwb$, hence pagewise. Since limits in $\lwbe$ are computed pagewise, we have $\Ker p\in\lwbe$. We use the characterization to exhibit a
counterexample to stability of fibrations under pullback along any map.

\begin{lem} Let $p\colon L\to M$ be a morphism in $\lwb$ with $L=(L_i,d^L_i,s^L_i,w^L_i)$ and $M=(M_i,d^M_i,s^M_i,w^M_i)$ in $\lwbs$. Then $\Ker p$ lies in $\lwbs$ if and only if for all $i\geq 0$  the following conditions hold
\begin{center}
\begin{tabular}{lllll}
$(C_\sigma^i) $& $\forall\ a\in L_i$, & $s^M_i p a=0$ & $\Longrightarrow$ & $\exists\ a'\in L_i,\; p a'=p a \text{ and }  s^L_i a'=0$, \\
$(C_\omega^{i+1})$ & $\forall\ a\in L_{i+1},$ & $w^M_{i+1} p a =0$ & $\Longrightarrow$ & $\exists\ a'\in L_{i+1},\; p a' =p a \text{ and }  w^L_{i+1} a'=0.$
\end{tabular}
\end{center}
\end{lem}

\begin{proof} We use the characterization from Proposition~\ref{P: charspseq} for a linear witness book $K\in\lwbe$ to be in $\lwbs$, that is: for all $i$,
$\Ker w^K_{i+1}=\Img s^K_i$ (which will 
 correspond to the condition $C_\sigma^{i}$) and $\Ker s^K_i=\Img w^K_{i+1}$ (which will correspond to condition $C_\omega^{i+1}$).
 
 We let $K=\Ker p$, and $w^K=w^L$,  $s^K=s^L$.

 Assume we have $\Ker w^K_{i+1}=\Img s^K_i$. In particular, for every $u\in L_{i+1}$, such that $pu=0$ and $w_{i+1}^Lu=0$, there exists 
 $v\in L_{i}$ such that $pv=0$ and $s_i^Lv=u$.
 Let $a\in L_i$ with $s^M_i p a=ps_i^L a=0$, and set $u=s_i^L a$. We  have $pu=0$ and $w_{i+1}^Lu=0$, hence there exists  
 $v\in L_{i}$ such that $pv=0$ and $s_i^Lv=u$. Let $a'=a-v$. Then $pa'=pa$ and $s_i^La'=s_i^La-s_i^Lv=u-u=0$. Hence condition $C_\sigma^i$ holds. 
 
 Conversely assume $C_\sigma^i$ holds. Let $u\in L_{i+1}$, such that $pu=0$ and $w_{i+1}^Lu=0$. Since $L\in\lwbs$,  there exists $v\in L_i$ such that $s_i^Lv=u$. In particular $s_i^Mpv=0$. Hence there exists $v'\in L_i$ with $pv'=pv$ and $s_i^L v'=0$. In particular $s_i(v-v')=u$ and $p(v-v')=0$.

The proof that condition $C_\omega^{i+1}$ is equivalent to $\Ker s^K_i=\Img w^K_{i+1}$  proceeds in exactly the same way, exchanging the roles of $s$ and $w$.
 \end{proof}
 
 \begin{prop} There is no model category structure on $\lwb$ or on $\lwbe$ such that the category of fibrant objects is the category $\lwbs$ with weak equivalences $\Ee_r$.
 \end{prop}
 
 \begin{proof} The proof is by contradiction. If such a model category structure exists with fibrant objects $\lwbs$, then any map $f\colon L\to M$ in $\lwb$ between two objects $L,M\in \lwbs$ factors as $f=pi$ with $i\colon L\to K$ an acyclic cofibration and $p\colon K\to M$ a fibration. Since $M$ is fibrant and $p$ is a fibration, $K$ is  fibrant.
 Since fibrations are stable under pullbacks, $p$ should satisfy conditions $C_\sigma^k$ and $C_\omega^k$ for every $k$.  
 
 We provide the following counterexample.
Let $f=\Nn(\pi)$ where $\pi\colon \Dd_r(0,0)\to \kk(0,0)$ is the projection in $\spse$. Recall that $\Nn\Dd_r(0,0)=\Yy(r,0,0)$ with generator in 
$\Yy(r,0,0)_r^{0,0}$ denoted $\alpha_r^{0,0}$. We have  for all $k$, $(\Nn\kk(0,0))_k=\kk^{0,0}$ with generator denoted $\beta_k^{0,0}$.
We have $s_k^{\Nn\kk(0,0)}=0$.
For every $0\leq k\leq r$, we have $\pi_k$ is surjective, hence $f_k$ is surjective by Proposition \ref{P:Nsurj}, and $f_r(\alpha_r^{0,0})=\beta_r^{0,0}$.  Finally $H(\Yy(r,0,0)_r)=0$.

Let us write $f=pi$, with $i\colon \Yy(r,0,0)\to K\in\Ee_r$, so that $H(K_r)=0$. Since $f_r$ is surjective, $p_r\colon K_r\to \kk^{0,0}$ is surjective and there exists $z\in K_r$ such that $p_rz=\beta_r^{0,0}$. In particular $s_rp_rz=0$. By condition $C^r_\sigma$,  there exists $z'\in K_r$ such that $p_rz'=\beta_r^{0,0}$ and $s_r^Kz'=0$.  But $\Ker s_r^K=\Ker d_r^K=\Img d_r^K$ because $K\in\lwbs\subset \lwbe$ and $H(K_r)=0$. Let $y\in K_r$ be such that $d_r^Ky=z'$. We have $p_rz'=\beta_r^{0,0}=d_rp_ry$, a contradiction.

Note that this proof by contradiction also works for $\lwbe$.
\end{proof}

\begin{cor} There is no model category structure on $\espse$ such that the category of fibrant objects is the category $\spse$ with weak equivalences $\Ee_r$.
\end{cor}

\begin{proof} This follows from the equivalences of categories with weak equivalences between $(\spse,\Ee_r)$ and $(\lwbs, \Ee_r)$ 
and between $(\espse,\Ee_r)$ and $(\lwbe, \Ee_r)$.
\end{proof}

\bibliographystyle{abbrv}

\bibliography{bibli_presheaves}

\end{document}